\documentclass[10pt, a4paper]{article}

\usepackage{graphicx}
\usepackage{color}
\usepackage{amsmath}
\usepackage{amssymb}
\usepackage{amscd}
\usepackage{bbm}
\usepackage{tikz}
\usetikzlibrary{patterns}
\usepackage{hyperref}
\hypersetup{
    bookmarks=true,         
    colorlinks=true,       
    linkcolor=red,          
    citecolor=green,        
    filecolor=magenta,      
    urlcolor=cyan           
}

\usepackage[utf8]{inputenc}
\usepackage{amsthm}
\usepackage{bbm} 
\usepackage[linesnumbered,lined,boxed,commentsnumbered]{algorithm2e}

\usepackage{marginnote}

\newtheorem{Theorem}{Theorem}
\newtheorem{Assumption}{Assumption}

\newtheorem{Lemma}{Lemma}
\newtheorem{Proposition}{Proposition}
\newtheorem{Corollary}{Corollary}
\newtheorem{Remark}{Remark}

\newcommand{\inr}[1]{\bigl< #1 \bigr>}

\newcommand{\norm}[1]{\left\|#1\right\|}%

\newcommand\eps{\epsilon}

\newcommand{\beginproof}{{\bf Proof. {\hspace{0.2cm}}}}
\def \endproof
{{\mbox{}\nolinebreak\hfill\rule{2mm}{2mm}\par\medbreak}}

\DeclareMathOperator*{\argmax}{argmax}

\def\ds1{\textrm{1\kern-0.25emI}} 

\newcommand \E{\mathbb{E}}

\newcommand \cA{{\cal A}}
\newcommand \cB{{\cal B}}

\newcommand \cE{{\cal E}}

\newcommand \cI{{\cal I}}

\newcommand \cK{{\cal K}}

\newcommand \cN{{\cal N}}
\newcommand \cO{{\cal O}}

\newcommand \cS{{\cal S}}

\newcommand \bE{{\mathbb E}}

\newcommand \bN{{\mathbb N}}

\newcommand \bP{{\mathbb P}}

\newcommand \bR{{\mathbb R}}

\usepackage[]{}

\DeclareMathOperator*{\Med}{Med}
\DeclareMathOperator*{\Tr}{Tr}

\DeclareMathOperator*{\Card}{Card}

\usepackage{braket}
\usepackage{authblk}

\begin{document}

\title{Robust subgaussian estimation of a mean vector in nearly linear time}
\author[1]{Jules Depersin and Guillaume Lecu{\'e}  \\ email: \href{mailto:jules.depersin@ensae.fr}{jules.depersin@ensae.fr}, email: \href{mailto:lecueguillaume@gmail.com}{guillaume.lecue@ensae.fr} \\ CREST, ENSAE, IPParis. 5, avenue Henry Le Chatelier, 91120 Palaiseau, France.}

\date{}                     
\setcounter{Maxaffil}{0}
\renewcommand\Affilfont{\itshape\small}

\maketitle

\begin{abstract}
We construct an algorithm,  running in time $\tilde\cO(N d + uK d)$ , which is robust to outliers and heavy-tailed data and which achieves the subgaussian rate from \cite{lugosi2019sub}
\begin{equation}\label{eq:intro_subgaus_rate}
\sqrt{\frac{\Tr(\Sigma)}{N}}+\sqrt{\frac{\norm{\Sigma}_{op}K}{N}}
\end{equation}with  probability at least $1-\exp(-c_0K)-\exp(-c_1 u)$ where $\Sigma$ is the covariance matrix of the \textit{informative data}, $K\in\{1, \ldots, K\}$ is some parameter (number of block means) and $u\in\bN^*$ is another parameter of the algorithm. This rate is achieved when $K\geq c_2 |\cO|$ where $|\cO|$ is the number of outliers in the database and under the only assumption that the informative data have a second moment. The algorithm is fully data-dependent and does not use in its construction the proportion of outliers nor the rate in \eqref{eq:intro_subgaus_rate}. Its construction combines recently developed tools for Median-of-Means estimators and covering-Semi-definite Programming \cite{MR3909640,PTZ12}. We also show that this algorithm can automatically adapt to the number of outliers.
\end{abstract}

\noindent\textbf{AMS subject classification:} 	62F35\\
\textbf{Keywords:} Robustness, algorithms, heavy-tailde data.

\section{Introduction on the robust mean vector estimation problem} 
\label{sec:introduction_on_the_mean_vector_problem}

Estimating the mean of a random variable in a $d$-dimensional space when given some of its realizations is arguably the oldest and most fundamental problem of statistics. In the past few years, it has received important attention from two communities: the Statistics \cite{AIHPB_2012__48_4_1148_0,minsker2015geometric,MR3845006,CG,lugosi2019sub,MR3851758,LMSL,hopkins2018sub,Bartlett19} and Computer Science  \cite{MR3631028,MR3945261,diakonikolas2018robustly,diakonikolas2016robust,MR3826316,MR3909639,MR3909640} communities. Both communities consider the problem of \textit{robust mean estimation}, focusing mainly on different definitions of robustness.

In recent years, many efforts have been made by the Statistics community on the construction of estimators performing in a \textit{subgaussian way} for heavy-tailed data. Such estimators achieve the same statistical properties as the empirical mean of a $N$-sample of i.i.d. gaussian variables $\cN(\mu, \Sigma)$ where $\mu\in\bR^d$  and $\Sigma\succeq0$ is the covariance matrix. In that case, for a given confidence $1-\delta$, the subgaussian rate as defined in \cite{lugosi2019sub} is  (up to an absolute multiplicative constant)
\begin{equation}\label{eq:sub_gauss}
 r_\delta = \sqrt{\frac{\Tr(\Sigma)}{N}} + \sqrt{\frac{||\Sigma||_{op} \log(1/\delta)}{N}}
\end{equation}where $\Tr(\Sigma)$ is the trace of $\Sigma$ and $||\Sigma||_{op}$ is the operator norm of $\Sigma$. Indeed, it follows from Borell-TIS's inequality (see Theorem~7.1 in \cite{Led01} or pages 56-57 in \cite{MR2814399}) that with probability at least $1-\delta$,
\begin{equation*}
\norm{\bar X_N-\mu}_2 = \sup_{\norm{v}_2\leq 1}\inr{\bar X_N-\mu,v}\leq \E \sup_{\norm{v}_2\leq 1}\inr{\bar X_N-\mu,v} + \sigma \sqrt{2\log(1/\delta)}
\end{equation*}where $\sigma =\sup_{\norm{v}_2\leq1} \sqrt{\E\inr{\bar X_N-\mu, v}^2}$. It is straightforward to check that $\E \sup_{\norm{v}_2\leq 1}\inr{\bar X_N-\mu,v}\leq \sqrt{\Tr(\Sigma)/N}$ and $\sigma=\sqrt{\norm{\Sigma}_{op}/N}$, which leads to the rate in \eqref{eq:sub_gauss} (up to the constant $\sqrt{2}$ on the second term in \eqref{eq:sub_gauss}). In most of the recent works, the effort has been made to achieve the rate $r_\delta$ for i.i.d. heavy-tailed data even under the minimal requirement that the data only have a second moment. Under this second-moment assumption only, the empirical mean cannot achieve the rate \eqref{eq:sub_gauss} and one needs to consider other procedures\footnote{Under only a second-moment assumption, the empirical mean  achieves the rate $\sqrt{\Tr(\Sigma)/(\delta N)}$ which can not be improved in general.}. Over the years, some procedures have been proposed to achieve such a goal: a Le Cam test estimator, called a tournament estimator in \cite{lugosi2019sub}, a minmax Median-Of-Means estimator in \cite{LMSL} and a PAC-Bayesian estimator in \cite{CG}. The first two one are based on the median-of-means principle that we will also use. 

On the other side, the Computer Science community mostly considers  a different definition of robustness and targets a different goal. In many recent CS papers, algorithms (not only estimators) have been constructed and proved to be robust with respect to a \textit{contamination} of the dataset that is when some of the data are replaced by other data which may have nothing to do with the original batch. This covers the Huber $\eps$-contamination model but also adversarial data which receives an important attention recently in the deep learning community. Moreover, the Computer Science community looks at the problem of robust mean estimation from  algorithmic perspectives such as the running time. A typical result in this line of research is Theorem~1.3 from \cite{MR3909640} that we recall now.

\begin{Theorem}[Theorem~1.3, \cite{MR3909640}]\label{theo:diakonikolas} Let $X_1, \ldots, X_N$ be random vectors in $\bR^d$. We assume that there is a partition $\{1, \ldots, N\}=\cO\cup\cI$ such that nothing is assumed on $(X_i)_{i\in\cO}$ and $(X_i)_{i\in\cI}$ are independent with mean $\mu$ and covariance matrix $\Sigma \preceq \sigma^2 I_d$. We assume that $\eps=|\cO|/N$ is such that $0<\eps<1/3$ and $N\gtrsim d \log(d) /\eps$. There exists an algorithm running in $\tilde \cO(Nd)/{\rm poly}(\eps)$ which outputs $\hat \mu_\eps$ such that with probability at least $9/10$, $\norm{\hat \mu_\eps - \mu}_2\lesssim \sigma \sqrt{\eps}$.
\end{Theorem} 

The first result proving the existence of a polynomial time algorithm robust to contamination may be found in \cite{MR3631028}. Theorem~\ref{theo:diakonikolas} improves upon many existing results since it achieves the optimal information theoretic-lower bound with a (nearly) linear-time algorithm.  

Finally, there are two recent papers for which both algorithmic and statistical considerations are important. In \cite{hopkins2018sub,Bartlett19}, algorithms achieving the subgaussian rate in \eqref{eq:sub_gauss} have been constructed. They both run in polynomial time : $\cO(N^{24} + Nd)$ for \cite{hopkins2018sub} and $\cO(N^4+N^2d)$ for \cite{Bartlett19} (see \cite{Bartlett19} for more details on these running times
). They do not consider a contamination of the dataset even though their results easily extend to this setup. Some other estimators which have been proposed in the Statistics literature are very fast to compute but they do not achieve the optimal subgaussian rate from \eqref{eq:sub_gauss}. A typical example is Minsker's geometric median estimator \cite{minsker2015geometric} which achieves the rate $\sqrt{\Tr(\Sigma) \log(1/\delta)/N}$ in linear time $\tilde \cO(Nd)$. All the later three papers use the Median-of-means principle. We will use this principle but only to construct a starting point (which will simply be the coordinate-wise median) and for the computation of the step size (where we will only use the one dimensional definition of the median along the descent line direction). What we mainly borrow from the literature on MOM estimators is the advantage to work with local block means instead of the data themselves. We will identify two such advantages by doing so: a stochastic one and a computational one (see Remark~\ref{rem:effects_blocks} below).  

Robust mean estimation have been raised in pioneered works in robust statistics from Huber \cite{MR0161415,MR2488795}, Tukey \cite{MR0120720,MR0133937} or Hampel \cite{MR0359096,MR0301858}. Their concerns was more about robustness to model misspecification and on the breakdown point property (``smallest amount of contamination necessary to upset an estimator entirely'' taken from \cite{MR1193313}). The computational problem connected to this issue was not of primary interest even though it was already raised, for instance, in Section~5.3 from \cite{MR1193313} for the construction of Tukey contours (a $d$-dimensional definition of  quantiles).

The aim of this work is to show that a single algorithm can answer the three problems: robustness to heavy-tailed data, to contamination and computational cost. In this article, we construct an algorithm running in time $\tilde\cO(N d + u\log(1/\delta) d)$ which outputs an estimator of the true mean achieving the subgaussian rate \eqref{eq:sub_gauss} with confidence $1-\delta$ (for $\exp(-c_0N)\leq \delta\leq \exp(-c_1|\cO|)$) on a corrupted database and under a second moment assumption only. It is therefore robust to heavy-tailed data and to contamination. Our approach takes ideas from both communities: the median-of-means principle which has been recently used in the Statistics community and a SDP relaxation from \cite{MR3909640} which can be computed fast. The baseline idea is to construct $K$ equal size groups of data from the $N$ given ones and to compute their empirical means $\bar X_k, k=1, \ldots, K$. These $K$ empirical means are  used successively to find a robust descent direction thanks to a SDP relaxation from  \cite{MR3909640}. We prove the robust subgaussian statistical property of the resulting  descent algorithm under the only following assumption.


\begin{Assumption}\label{assum:first}There exists a partition $\cI\cup\cO=\{1, \ldots, N\}$ of the dataset $(X_i)_{i \leq N}$ such that 1) nothing is assumed on $(X_i)_{i\in\cI}$ 2) $(X_i)_{i\in\cI}$ are independent with mean $\mu$ and covariance $\E (X_i-\mu) (X_i-\mu)^\top) \preceq \Sigma $ where $\Sigma$ is a given (unknown) covariance matrix.
\end{Assumption}

Assumption~\ref{assum:first} covers the two concepts of robustness considered in the Statistics and Computer Science communities since the \textit{informative data} (data indexed by $\cI$) are only assumed to have a second moment and there are $|\cO|$ outliers onto which we do not make any assumption. Our aim is to show that the rate of convergence \eqref{eq:sub_gauss} which is the rate achieved by the empirical mean in the ideal i.i.d. Gaussian case can be achieved in the corrupted and heavy-tailed setup from Assumption~\ref{assum:first} with a fast algorithm.\\

The paper is organized as follows. In the next section, we give a high-level description of the algorithm and its statistical and computation performances. In section~3, we prove its statistical properties and give a precise definition of the algorithm. In Section~4, we study the statistical performance of the SDP relaxation at the heart of the descent direction. In Section~5, we fully characterize its computational cost. In Section~\ref{sec:adaptive_choice_of_}, we construct a procedure achieving the same statistical properties and can automatically adapt to the number of outliers.

 
 
 
\section{Construction of the algorithms and main result} 
\label{sec:construction_of_the_algorithms_and_main_result}
The construction of our robust subgaussian descent procedure is using two ideas. The first one comes from the median-of-means (MOM) approach which has recently received a lot of attention in the statistical and machine learning communities  \cite{MR3124669,LO,MR3576558,MS,minsker2015geometric}. The MOM  approach \cite{MR702836,MR1688610,MR855970,MR762855} often yields robust estimation strategies (but usually at a high computational cost). Let us give the general idea behind that approach: we first randomly split the data into $K$ equal-size blocks $B_1,\ldots ,B_K$ (if $K$ does not divide $N$, we just remove some data). We then compute the empirical mean within each block: for $k=1,\ldots,K$,
\begin{equation*}
 \bar{X}_k=\frac{1}{|B_k|}\sum_{i \in B_k} X_i
 \end{equation*}  
 where we set $|B_k|=\Card(B_k)=N/K$. In the one-dimensional case, we then take the median of the latter $K$ empirical means to construct a robust \emph{and subgaussian} estimator of the mean \cite{MR3576558}. It is more complicated in the multi-dimensional case, where there is no \textit{definitive} equivalent of the one dimensional median but several candidates:  coordinate-wise median, the geometric median (also known as Fermat point), the Tukey Median,  among many others (see \cite{small1990survey}). The strength of this approach is the robustness of the median operator, which leads to good statistical properties even on corrupted databases.  For the construction of our algorithm, we actually only use  the idea of grouping the data and computing their $K$ means $\bar X_k, k=1, \ldots, K$. 

 Finding good descent directions in the heavy-tailed and corrupted scenario considered in Assumption~\ref{assum:first} in reasonnable time is a main issue. A construction has been proposed by  \cite{Bartlett19} which also uses a SDP relaxation, which costs $\cO(N^4+Nd)$ to be computed. Our approach also uses a SDP relaxation, with an other SDP. It is based on the observation that $\mu$ is solution of the minimization problem $\min_{\nu\in\bR^d}f(\nu)$ where $f:\nu\in\bR^d\to \norm{\E X - \nu}_2^2$ and $X$ is any random vector with mean $\mu$. One way to approach $\mu$ is therefore to run a gradient descent algorithm using $f$ as an objective function:  from $x_c\in\bR^d$ we go to the next iteration with $x_c- \theta \nabla f(x_c)$ where $\theta\geq0$ is a step size. Since $\nabla f(x_c) = x_c-\bE X$, for $\theta=1$, the latter algorithm achieves the target mean $\mu$ in one step, which is not surprising given that $x_c-\E X$ is the best descent direction towards $\E X$ starting from $x_c$. We can also re-write that as a matrix problem : the top eigenvector of
 \begin{equation}\label{eq:minimizing_pb}
 \argmax_{M\succeq 0, \Tr(M)=1}\inr{M, (\E X-x_c)(\E X-x_c)^\top}
 \end{equation}is given by $\frac{x_c-\E X}{\norm{x_c-\E X}_2}$, which is the best descent direction we are looking for. \\ 

  Of course, we don't know $(\E X-x_c)(\E X-x_c)^\top$ in \eqref{eq:minimizing_pb} but we are given a database of $N$ data $X_1, \ldots, X_N$ (among which $|\cI|$ of them have mean $\mu$). We use these data to estimate in a robust way the unknown quantity $(\E X-x_c)(\E X-x_c)^\top$ in \eqref{eq:minimizing_pb}. Ideally, we would like to identify the \textit{informative data} and then use $(1/|\cI|)\sum_{i\in\cI}(X_i-x_c)(X_i-x_c)^\top$ or its block means version $(1/|\cK|)\sum_{k\in\cK}(\bar X_k-x_c)(\bar X_k-x_c)^\top$, where $\cK=\{k:B_k\cap \cO=\emptyset\}$, to estimate this quantity   but this information is not available either.  

  To address this problem we use a tool introduced in \cite{MR3909640} adapted to the block means. The idea is to endow each block mean $\bar{X}_k$ with a weight $\omega_k$ taken in $\Delta_{K}$ defined as
 \begin{equation*}
 \Delta_{K}=\left\{( \omega_k)_{k=1}^K: 0 \leq \omega_k \leq \frac{1}{9K/10}, \sum_{k=1}^K \omega_k =1   \right\}.
 \end{equation*}Ideally we would like to put $0$ weights to all block means $\bar{X}_k$ corrupted by an outliers. But, we cannot do it since $\cK$ is unknown. To overcome this issue, we learn the optimal weights and consider the following minmax optimization problem
 \begin{equation}\label{formule1} \tag{$E_{x_c}$}
 \underset{M \succeq 0, \Tr(M)=1}{\text{max}} \ \underset{w \in \Delta_{K}}{\text{min}} \ \inr{M, \sum_{k=1}^K \omega_k (\bar{X}_k-x_c)(\bar{X}_k-x_c)^\top}.
\end{equation}This is the dual problem from \cite{MR3909640} adapted to the block means. The key insight from \cite{MR3909640} is that  an approximating solution $M_c$ of the maximization problem in \eqref{formule1} can be  obtained in reasonable time using a covering SDP approach \cite{MR3909640,PTZ12} (see Section~\ref{sec:SDP}). We expect a solution (in $M$) to \eqref{formule1} to be close to a solution of the minimization problem in \eqref{eq:minimizing_pb} -- which is $M^*=(\mu-\nu)(\mu-\nu)^\top/\norm{\mu-\nu}_2^2$ -- and the same for their top eigenvectors (up to the sign). 

At a high level description, the robust descent algorithm we perform outputs $\hat \mu_K $ after at most $\log d$ iterations of the form $x_c - \theta_c v_1$ where $v_1$ is a top eigenvector of an approximating solution $M_c$ to the problem \eqref{formule1} and $\theta_c$ is a step size. It starts at the coordinate-wise median of the means $\bar X_1, \ldots, \bar X_K$ . In Algorithm~\ref{algo:final}, we define precisely  the step size and the stopping criteria we use to define the algorithm (it requires too many notation to be defined at this stage). This algorithm outputs the vector $\hat \mu_K$ : its running time and statistical performances are gathered in the following result.

\begin{Theorem}\label{theo:main}
Grant Assumption~\ref{assum:first}. Let $K\in\{1,\ldots, N \}$ be the number of equal-size blocks and assume that $K\geq 300 |\cO|$. Let $u\in\bN^*$ be a parameter of the covering SDP used at each descent step. With probability at least $1-\exp(-K/180000)-(1/10)^u$, the descent algorithm finishes in $\tilde \cO(Nd+K u d)$ and outputs $\hat \mu_K$ such that 
\begin{equation*}
\norm{\hat \mu_K - \mu}_2\leq 808 \left(1200 \sqrt{\frac{\Tr(\Sigma)}{N}} +  \sqrt{\frac{1200\norm{\Sigma}_{op}K}{N}}\right).
\end{equation*}
\end{Theorem}

To make the presentation of the proof of Theorem~\ref{theo:main} as simple as possible we did not optimize the constants. Theorem~\ref{theo:main} generalizes and improves Theorem~\ref{theo:diakonikolas} in several ways. We first improve the confidence from a constant ``$9/10$'' to an exponentially large confidence $1-\exp(-c_0K)$. We obtain the result for any covariance structure $\Sigma$ and $\hat \mu_K$ does not require the knowledge of $\Sigma$ for its construction. We obtain a result which holds for any $N$ (even under the sample complexity). The construction of $\hat \mu_K$ does not require the knowledge of the exact proportion of outliers $\eps$ in the dataset unlike $\hat\mu_\eps$ in Theorem~\ref{theo:diakonikolas}. We only need to know that $K\gtrsim |\cO|$.  Moreover, using a Lepskii adaptation method it is also possible to automatically choose $K$ and therefore to adapt to the proportion of outliers if we have some extra knowledge on $\Tr(\Sigma)$ and $\norm{\Sigma}_{op}$ (see Section~\ref{sec:adaptive_choice_of_} for more details). Moreover, if we only care about constant $9/10$  confidence, our runtime does not depend on $\epsilon$ and is nearly-linear $\tilde \cO(Nd)$. We also refer the reader to Corollary~\ref{coro:lepski} for more comparison with Theorem~\ref{theo:diakonikolas}.

\begin{Remark}[Nearly-linear time]
We identify two important situations where the algorithm from Theorem~\ref{theo:main} runs in nearly-linear time that is in $\tilde\cO(Nd)$. First, when the number of outliers is known to be less than $\sqrt{N}$, we can choose $K\leq \sqrt{N}$ and $u=K$. In that case, the algorithm runs in $\tilde\cO(Nd)$ and the subgaussian rate is achieved with probability at least  $1-2\exp(-c_0K)$ for some constant $c_0$ (see also Corollary~\ref{coro:lepski_2_subgaussian} for an adaptive to $K$ version of this result). Another widely investigated situation is when we only want to have a constant confidence like $9/10$. In that case, one may chose $u=1$ and any values of $K\in[N]$ can be chosen (so we can have any number of outliers) to achieve the subgaussian rate with constant probability and in nearly-linear time $\tilde\cO(Nd)$ (see also Corollary~\ref{coro:lepski} for an adaptive to $K$ version of this result). 
\end{Remark}

Theorem~\ref{theo:main} improves the result from \cite{hopkins2018sub,Bartlett19} since $\hat \mu_K$ runs faster than the polynomial times $\cO(N^{24} + Nd)$ and  $\cO(N^4 + Nd)$ in \cite{hopkins2018sub} and \cite{Bartlett19}. The algorithm  $\hat \mu_K$ also does not require the knowledge of $\Tr(\Sigma)$ and $\norm{\Sigma}_{op}$. Finally, Theorem~\ref{theo:main} provides running time guarantees on the algorithm unlike in \cite{lugosi2019sub,LMSL,CG} and it improves upon the statistical performances from \cite{minsker2015geometric}.

\section{Proof of the statistical performance in Theorem~\ref{theo:main}} 
\label{sec:proof_of_the_statistical_performance_in_theorem_theo:main}
In this section, we prove the statistical performance of $\hat \mu_K$ as stated in Theorem~\ref{theo:main}. We first identify an event $\cE$ onto which we will derive the rate of convergence of the order of \eqref{eq:sub_gauss}. This event is also used to compute the running time of $\hat \mu_K$ in the next section as announced in Theorem~\ref{theo:main}.

\begin{Proposition}\label{Prop:MatriceLemme}
Denote by $\cE$ the event onto which for all matrix $M \succeq 0$ such that $\Tr(M)=1$, there are at least $9K/10$ of the blocks for which $\norm{M^{1/2} (\bar{X}_k- \mu)}_2 \leq  8r$ where 
\begin{equation}\label{eq:def_r}
 r = 1200 \sqrt{\frac{\Tr(\Sigma)}{N}} + \sqrt{\frac{1200\norm{\Sigma}_{op}K}{N}}.
 \end{equation} If Assumptions~\ref{assum:first} holds and  $K\geq 300|\cO|$ then  $\bP[\cE]\geq 1-\exp(-K/180000)$.
\end{Proposition}

Proposition~\ref{Prop:MatriceLemme} contains all the stochastic arguments we will use in this paper (constants have not been optimized). In other words, after identifying $\cE$ all the remaining arguments do not involve any other stochastic tools. Before proving Proposition~\ref{Prop:MatriceLemme}, let us first state a result that is of particular interest beyond our problem.

\begin{Corollary}\label{coro:dual_isometry}
On the event $\cE$, for all $M\in\bR^{d\times d}$ such that $M\succeq 0$ and $\Tr(\Sigma)=1$ there are at least $9K/10$ blocks such that for all $x_c\in\bR^d$,
\begin{equation}\label{eq:coro_dual_isometry}
\norm{M^{1/2}(\mu - x_c)}_2 - 8r\leq \norm{M^{1/2}(\bar X_k - x_c)}_2\leq \norm{M^{1/2}(\mu - x_c)}_2 + 8r.
\end{equation}
\end{Corollary}

Let us now turn to a proof of Proposition~\ref{Prop:MatriceLemme}.  We first remark that if we were to only consider matrices $M$ of rank $1$, Proposition~\ref{Prop:MatriceLemme} would boil down to show that for all $v\in\cS_2^{d-1}$ (the unit sphere in $\ell_2^d$) on more than $9/10$ blocks $|\inr{v,\bar X_k-\mu}|\leq 8r$. This is a ``classical'' result in the MOM literature which has been proved in \cite{lugosi2019sub} and \cite{LMSL}. We recall now this result and the short proof from  \cite{LMSL} for completeness. We will use it to prove Proposition~\ref{Prop:MatriceLemme}.


\begin{Lemma}\label{lemm:VecteurLemme}
Grant Assumption~\ref{assum:first} and assume that $K\geq 300|\cO|$. With probability at least $1- \exp(-K/180000)$, for all $v\in\cS_2^{d-1}$, there are at least $99K/100$ of the blocks $k$ such that $|\inr{v,\bar{X}_k-\mu}| \leq r$.
\end{Lemma}
\beginproof
We want to show that with probability at least $1-\exp(-K/180000)$, for all $v\in \cS_2^{d-1}$, 
 \begin{equation*}
 \sum_{k\in[K]} I(|\inr{\bar{X}_k - \mu, v}|> r)\leq K/100.
 \end{equation*}We take $\cK=\{k\in[K]: B_k\cap \cO=\emptyset\}$. We define $\phi(t) = 0 $ if $t\leq1/2$, $\phi(t) = 2(t-1/2)$ if $1/2\leq t\leq 1$ and $\phi(t) = 1$ if $t\geq1$. We have $I(t\geq1)\leq \phi(t)\leq I(t\geq1/2)$ for all $t\in\bR$ and so
  \begin{align*}
 &\sum_{k\in \cK} I(|\inr{\bar{X}_k - \mu, v}|> r)\leq \sum_{k\in\cK} I(|\inr{\bar{X}_k - \mu, v}|> r) - \bP[|\inr{\bar{X}_k - \mu, v}|> r/2] + \bP[|\inr{\bar{X}_k - \mu, v}|> r/2]\\
 &\leq \sum_{k\in\cK}\phi\left(\frac{|\inr{\bar{X}_k - \mu, v}|}{r}\right) - \bE \phi\left(\frac{|\inr{\bar{X}_k - \mu, v}|}{r}\right) + \bP[|\inr{\bar{X}_k - \mu, v}|> r/2]\\
& \leq \sup_{v\in \cS_2^{d-1}}\left(\sum_{k\in\cK}\phi\left(\frac{|\inr{\bar{X}_k - \mu, v}|}{r}\right) - \bE \phi\left(\frac{|\inr{\bar{X}_k - \mu, v}|}{r}\right) \right)+  \sum_{k\in\cK} \bP[|\inr{\bar{X}_k - \mu, v}|> r/2].
 \end{align*}For all $k\in\cK$, we have
 \begin{align*}
 \bP[|\inr{\bar{X}_k - \mu, v}|> r/2]\leq \frac{\bE \inr{\bar{X}_k - \mu, v}^2}{(r/2)^2} \leq \frac{4Kv^\top \Sigma v}{Nr^2}\leq \frac{4K\sup_{v\in \cS_2^{d-1}}v^\top \Sigma v}{Nr^2} = \frac{4K\norm{\Sigma}_{op}}{Nr^2}\leq \frac{1}{300}
 \end{align*}because $r^2\geq 1200K \norm{\Sigma}_{op}/N$. Next, using the bounded difference inequality (Theorem~6.2 in \cite{MR3185193}), the symmetrization argument and the contraction principle (Chapter~4 in \cite{MR2814399}), with probability at least $1-\exp(-K/180000)$, 
 \begin{align*}
  &\sup_{v\in \cS_2^{d-1}}\left(\sum_{k\in\cK}\phi\left(\frac{|\inr{\bar{X}_k - \mu, v}|}{r}\right) - \bE \phi\left(\frac{|\inr{\bar{X}_k - \mu, v}|}{r}\right) \right)\\
  & \leq \bE \sup_{v\in S}\left(\sum_{k\in\cK}\phi\left(\frac{|\inr{\bar{X}_k - \mu, v}|}{r}\right) - \bE \phi\left(\frac{|\inr{\bar{X}_k - \mu, v}|}{r}\right) \right)+ \sqrt{\frac{|\cK|K}{360000}}\\
  &\leq \frac{4K}{Nr} \bE \sup_{v\in \cS_2^{d-1}} \inr{v, \sum_{i\in \cup_{k\in\cK}B_k}\eps_i (X_i-\mu)} +  \sqrt{\frac{|\cK|K}{360000}}\\
  & =  \frac{4K}{\sqrt{N}r}\bE \norm{\frac{1}{\sqrt{N}}\sum_{i\in \cup_{k\in\cK}B_k}\eps_i (X_i-\mu)}_{2}  +  \sqrt{|\cK|K/360000}\leq \frac{K}{300}
 \end{align*}because $r\geq 1200 \bE \norm{\sum_{i\in \cup_{k\in\cK}B_k}\eps_i (X_i-\mu^*)}_2/\sqrt{N}$ since
 \begin{equation*}
 \bE \norm{\frac{1}{\sqrt{N}}\sum_{i\in \cup_{k\in\cK}B_k}\eps_i (X_i-\mu)}_{2}\leq \sqrt{\bE \norm{\frac{1}{\sqrt{N}}\sum_{i\in \cup_{k\in\cK}B_k}\eps_i (X_i-\mu)}_2^2}=\sqrt{\frac{|\cup_{k\in\cK}B_k|}{N}}\sqrt{\Tr(\Sigma)}\leq \sqrt{\Tr(\Sigma)}.
 \end{equation*}

 As a consequence, when $K\geq 300 |\cO|$, with probability at least $1-\exp(-K/180000)$, for all $v\in \cS_2^{d-1}$, 
 \begin{equation*}
 \sum_{k\in [K]} I(|\inr{\bar{X}_k - \mu, v}|> r)\leq |\cO| + \frac{|\cK|}{300} + \frac{K}{300}\leq \frac{K}{100}.
 \end{equation*}

\endproof


\textbf{Proof of Proposition~\ref{Prop:MatriceLemme}:} Let $M\in\bR^{d\times d}$ be such that $M\succeq 0$ and $\Tr(\Sigma)=1$. Denote by $\mathcal{A}_M=\{k \in [K]: \norm{M^{1/2}(\bar{X}_k-\mu)}_2 \geq  8 r\}$ and assume that $|\cA_M|\geq 0.1K$.  Let $G$ be a Gaussian vector in $\bR^d$ with mean $0$ and covariance matrix $M$ (and independent from $X_1, \ldots, X_N$). We consider the random variable $Z=\sum_{k\in[K]} I\left(|\inr{\bar{X}_k-\mu,G}| > 5r\right)$. We work conditionally to $X_1, \ldots, X_N$ in this paragraph. For all $k\in[K]$, $ \inr{\bar{X}_k-\mu,G} $ is a centered Gaussian variable with variance $\sigma_k^2:=\norm{M^{1/2}(\bar{X}_k-\mu)}_2^2$. In particular, for all $k\in\cA_M$, if we denote by $g$ a standard real-valued Gaussian variable,  we have $\bP_G\left[|\inr{\bar{X}_k-\mu,G}| >5r\right]\geq \bP_G\left[|\inr{\bar{X}_k-\mu,G}| > 5\sigma_k/8\right] = 2 \bP[g>5/8]\geq 0.528$ (where $\bP_G$ (resp. $\bE_G$) denotes the probability (resp. expectation) w.r.t. $G$ conditionally on $X_1, \ldots, X_N$). Hence, $\E_G Z\geq 0.528|\cA_M|\geq 0.0528 K$. Since $|Z|\leq K$ a.s., it follows from Paley-Zygmund inequality (see Proposition~3.3.1 in \cite{MR1666908}) that 
\begin{equation*}
 \bP_G[Z> 0.01K]\geq \frac{(\E_G Z-0.01K)^2}{\E_G Z^2}\geq (0.0428)^2 =0.0018.
 \end{equation*} Moreover, it follows from the Borell-TIS inequality (see Theorem~7.1 in \cite{Led01} or pages 56-57 in \cite{MR2814399}) that with probability at least $1-\exp(-8)$, $\norm{G}_2\leq \E \norm{G}_2 + 4\sqrt{\norm{M}_{op}}$. Moreover, $\E\norm{G}_2 \leq  \sqrt{\Tr(M)}\leq 1$ and $\norm{M}_{op} \leq\Tr(M)\leq1$, so $\norm{G}_2\leq 5$ with probability at least  $1- \exp(-8) \geq 0.9996$. Since $0.9996+0.0018 > 1$ there exists a vector $G_M\in\bR^d$ such that $\norm{G_M}_2\leq 5$ and $\sum_{k\in[K]} I\left(|\inr{\bar{X}_k-\mu,G_M}| > 5 r\right) > 0.01 K$. We recall that this latter result holds when we assume that $|\cA_M|\geq 0.1K$.

Next, we denote by $\Omega_0$ the event onto which for all $v\in\cS_2^{d-1}$, there are at least $99K/100$ blocks such that $|\inr{\bar X_k-\mu,v}|\leq r$. We know from Lemma~\ref{lemm:VecteurLemme} that $\bP[\Omega_0]\geq 1-\exp(-K/180000)$. Let us place ourselves on the event $\Omega_0$ up to the end of the proof. Let $M\in\bR^{d\times d}$ be such that $M\succeq 0$ and $\Tr(\Sigma)=1$ and assume that $|\cA_M|\geq 0.1K$. It follows from the first paragraph of the proof that there exists $G_M\in\bR^d$ such that $\norm{G_M}_2\leq 5$ and $\sum_{k\in[K]} I\left(|\inr{\bar{X}_k-\mu,G_M}| > 5 r\right) >0.01 K$. Given that we work on the event $\Omega_0$, we have for  $v_M=G_M/\norm{G_M}_2$, that for more than $99K/100$ blocks  $|\inr{\bar{X}_k-\mu,v_M}| \leq r$ and so $|\inr{\bar{X}_k-\mu,G_M}| \leq \norm{G_M}_2r\leq 5r$ which contradicts the fact that $\sum_{k\in[K]} I\left(|\inr{\bar{X}_k-\mu,G_M}| > 5 r\right) >0.01 K$. Therefore, we necessarily have $|\cA_M|\leq 0.1K$, which concludes the proof.
 \endproof

 \textbf{Proof of Corollary~\ref{coro:dual_isometry}:}
Let us assume that the event $\cE$ holds up to the end of the proof. Let $M\in\bR^{d\times d}$ be such that $M\succeq0$ and $\Tr(\Sigma)=1$. Let $\cK_M = \{k\in[K]: \norm{M^{1/2}(\bar X_k - \mu)}_2\leq 8r\}$. On the event $\cE$, we have $|\cK_M|\geq 9K/10$.  Let $x_c\in\bR^d$. For all $k\in\cK_M$, we have $\norm{M^{1/2}(\mu-x_c)}_2\leq 8r$ and so
\begin{align*}
\norm{M^{1/2}(\bar X_k - x_c)}_2 &\in\left[\norm{M^{1/2}(\bar X_k - \mu)}_2 - \norm{M^{1/2}(\mu-x_c)}_2, \norm{M^{1/2}(\bar X_k - \mu)}_2 + \norm{M^{1/2}(\mu-x_c)}_2\right]\\
&\subset \left[\norm{M^{1/2}(\bar X_k - \mu)}_2 - 8r, \norm{M^{1/2}(\bar X_k - \mu)}_2 + 8r\right].
\end{align*}
 \endproof
 
Let us now turn to the study of the optimization problem \eqref{formule1} on the event $\cE$. Like in \cite{MR3909640}, we denote by  $OPT_{x_c}$ the optimal value of \eqref{formule1} and by $h_{x_c}: M \to  \underset{w \in \Delta_{K}}{\text{min}} \ \braket{M, \sum_k \omega_k (\bar{X}_k-x_c)(\bar{X}_k-x_c)^\top}$ its objective function to be minimized over the constraint set $\{M\in\bR^{d\times d}:M\succeq 0, \Tr(M)=1\}$.

\begin{Remark}\label{rem:objective_val}
For a given $M$, the optimal choice of $w\in\Delta_K$ in the definition of $h_{x_c}(M)$ is straightforward: one just have to put the maximum possible weight on the $9K/10$ smallest $\inr{M,  (\bar{X}_k-x_c)(\bar{X}_k-x_c)^\top}, k\in[K]$. Formally, we set $\mathcal{S}_M= \sigma(\{1,2,\cdots, 9K/10 \})$, where $\sigma$ is a permutation on $[K]$ that arranges the $(\bar{X}_k- x_c)^\top M  (\bar{X}_k- x_c), k\in[K]$  in ascending order: 
\begin{align*}
(\bar{X}_{\sigma(1)}- x_c)^\top M  (\bar{X}_{\sigma(1)}- x_c) \leq (\bar{X}_{\sigma(2)}- x_c)^\top & M  (\bar{X}_{\sigma(2)}- x_c) \leq \cdots\leq (\bar{X}_{\sigma(K)}- x_c)^\top M  (\bar{X}_{\sigma(K)}- x_c).
\end{align*}
Then we get $h_{x_c}(M)=(1/|\mathcal{S}_M|)\sum_{k \in \mathcal{S}_M} (\bar{X}_k- x_c)^\top M  (\bar{X}_k- x_c)$. 
\end{Remark}


The first lemma deals with the optimal value of \eqref{formule1} when the current point $x_c$ is far from $\mu$.

\begin{Lemma}\label{DistanceLemme}
On the event $\cal E$, for all $x_c\in\bR^d$, if $\norm{x_c-\mu}_2> 16r$ then $$(8/9)(\norm{x_c-\mu}_2-8r)^2\leq OPT_{x_c}  \leq (\norm{x_c-\mu}_2+8r)^2.$$
\end{Lemma}

\beginproof
 Let $M$ be a matrix such that $M \succeq 0$ and  $\Tr(M)=1$. Set $\cK_M = \{k\in[K]: \norm{M^{1/2}(\bar X_k - \mu)}_2\leq 8 r\}$. On the event $\cE$, we have $|\cK_M|\geq 9K/10$ and it follows from the proof of Corollary~\ref{coro:dual_isometry} that for all $k\in\cK_M$ and all $x_c\in\bR^d$, 
 \begin{equation}\label{eq:inter_lemma_distance_opt}
\norm{M^{1/2}(\mu - x_c)}_2 - 8r\leq \norm{M^{1/2}(\bar X_k - x_c)}_2\leq \norm{M^{1/2}(\mu - x_c)}_2 + 8r.
\end{equation} Then we define a weight vector $ \tilde \omega \in \Delta_{K}$ by setting for all $k\in[K]$
\begin{equation*}
\tilde \omega_k = \left\{
    \begin{array}{ll}
        1/|\cK_M|& \mbox{if } k \in \cK_M\\
        0 & \mbox{else.}
    \end{array}
\right.
\end{equation*}It follows from the definition of $h_{x_c}$ and \eqref{eq:inter_lemma_distance_opt} that
\begin{equation}\label{eq:inter_2_lemma_distance_opt}
h_{x_c}(M) \leq  \sum_{k\in[K]} \tilde \omega_k (\bar{X}_k- x_c)^\top M  (\bar{X}_k- x_c) = \frac{1}{|\cK_M|}\sum_{k\in\cK_M}\norm{M^{1/2}(\bar X_k-x_c)}_2^2 \leq \left(\norm{M^{1/2}(\mu - x_c)}_2 + 8r\right)^2.
\end{equation}Taking the maximum over all $M\in\bR^d$ such that $M\succeq 0$ and $\Tr(\Sigma)=1$ on both side of  the latter inequality yields the right-hand side inequality of Lemma~\ref{DistanceLemme}.

For the left-hand side inequality of Lemma~\ref{DistanceLemme}, we let $x_c\in\bR^d$ be such that $\norm{x_c-\mu}_2> 16r$. Let  $M$ be such that $M\succeq 0$ and $\Tr(M)=1$.  We use the notation and observation from Remark~\ref{rem:objective_val}: we note that $|\cK_M \cap \cS_M| \geq 8K/10$ so that it follows from Corollary~\ref{coro:dual_isometry} that 
 \begin{align*}
 h_{x_c}(M)&=\frac{1}{9K/10}\sum_{k \in \mathcal{S}_M} \norm{M^{1/2} (\bar{X}_k- x_c)}_2^2 \geq \frac{1}{9K/10} \sum_{k \in \mathcal{A}_M \cap \mathcal{S}_M} \norm{ M^{1/2} (\bar{X}_k- x_c)}_2^2\\
  &\geq \frac{8K/10}{9K/10}  \left( \norm{M^{1/2}(\mu- x_c)}_2- 8r\right)^2.
 \end{align*} Then, taking the maximum over all  $M\succeq0$ such that $\Tr(M)=1$ on both sides, finishes the proof. 
\endproof

Next lemma shows that the top eigenvector of an approximating solution to \eqref{formule1} is aligned with the best possible descent direction $(\mu-x_c)/\norm{\mu-x_c}_2$. It is taken from the proof of Lemma~3.3 in \cite{MR3909640}. We reproduce here a short proof for completeness.

\begin{Proposition}\label{Prop:direction}
On the event $\cE$, if $M$ is a matrix such that $M \succeq 0$, $\Tr(M)=1$  and $h_{x_c}(M) \geq (\beta \norm{x_c-\mu}_2 +8r)^2$ for some $1/\sqrt{2}\leq \beta\leq1$, then any top eigenvector $v_1$ of $M$ satisfies
$$ \left|\inr{v_1,\frac{x_c-\mu}{\norm{x_c-\mu}_2}} \right| >\sqrt{2\beta^2-1}.$$
\end{Proposition}

\beginproof
 Let $M$ be a matrix such that $M \succeq 0$ , $\Tr(M)=1$ and $h_{x_c}(M) \geq (\beta \norm{x_c-\mu}_2 +8r)^2$ for some $1/\sqrt{2}\leq \beta\leq1$. We know from the proof of Lemma~\ref{DistanceLemme} (see Equation~\eqref{eq:inter_2_lemma_distance_opt}) that
$ h_{x_c}(M) \leq \left(\norm{M^{1/2} (\mu- x_c)}_2+8r \right)^2$. This implies that $\norm{M^{1/2}(\mu- x_c)}_2^2 \geq  \beta^2 \norm{\mu-x_c}_2^2$.

Let $\lambda_1 \geq \lambda_2 \geq \ldots \geq \lambda_d\geq 0$ denote the eigenvalues of $M$ and let $v_1,\ldots,v_d$ denote  corresponding eigenvectors. The conditions on $M$ implies that $\sum_{j} \lambda_j = 1$ and $\cB_M=(v_1, \ldots, v_d)$ is an orthonormal basis of $\bR^d$. We denote $v=(\mu-x_c) / \norm{\mu-x_c}_2$. We decompose $v$ in $\cB_M$ as $v= \sum_j \alpha_j v_j$  with $\sum_j \alpha_j^2 = 1$. Using this decomposition, we have $ v^\top M v = \sum_j \lambda_j \alpha_j^2$. We have $\lambda_1 = \lambda_1  \sum_j \alpha_j^2 \geq \sum_j \lambda_j \alpha_j^2 \geq  \beta^2$, so  $\lambda_1 \geq \beta^2$. Moreover, since $\sum_{j} \lambda_j = 1$, we have $\beta^2 \sum_j \alpha_j^2 \leq \sum_j \lambda_j \alpha_j^2 \leq \lambda_1 \alpha_1^2 +(1-\lambda_1)(1-\alpha_1^2) \leq \alpha_1^2+ (1-\beta^2)\sum_j \alpha_j^2$, so we have $\alpha_1^2 \geq (2\beta^2-1)$. As we know that $\alpha_1=\inr{v_1, v}$, we get the result. 
\endproof

Proposition~\ref{Prop:direction} is the first tool we need to construct a descent algorithm since it provides a descent/ascent direction (depending on the sign of the top eigenvector of an approximate solution to \eqref{formule1}). It remains to specify three other quantities to fully characterize our algorithm: a starting point, a step size and a stopping criteria. We start with the starting point. Here we simply use the coordinate-wise median-of-means. The following statistical guarantee on the coordinate-wise median-of-means is known or folklore but we want to put forward that in our case it holds on the event $\cE$. This again shows that $\cE$ is the only event we need to fully analyze all the building blocks of our algorithm. We recall that the coordinate-wise median-of-means is the estimator $\hat\mu^{(0)}\in\bR^d$ whose coordinates are for all $j\in[d], \hat \mu^{(0)}_j = {\rm med}(\bar X_{k,j}:k\in[K])$ where $\bar X_{k,j}$ is the $j$-th coordinate of the block mean $\bar X_k$ for all $k\in[K]$.

\begin{Proposition}\label{prop:coordinate_wise_MOM}
On the even $\cE$, we have $\norm{\hat \mu^{(0)} - \mu}_2\leq 8\sqrt{d} r$.
\end{Proposition}
\beginproof
Let us place ourselves on the event $\cE$ during all the proof. For all direction, $v\in\cS_2^{d-1}$, there are at least $9K/10$ blocks $k$ such that $|\inr{\bar X_k-\mu, v}|\leq 8 r$. In particular, for all $j\in[d], |\inr{\bar X_k - \mu, e_j}|\leq 8r$ where $(e_1, \ldots, e_d)$ is the canonical basis of $\bR^d$. That is for at least $9K/10$ blocks $|\bar X_{k,j}-\mu_j|\leq 8r$. In particular, the latter result is true for the median of $\{\bar X_{k,j}:k\in[K]\}$ that is for $\hat \mu^{(0)}_j$. We therefore have $\norm{\hat \mu^{(0)} - \mu}_\infty\leq 8r$ and so $\norm{\hat \mu^{(0)} - \mu}_2\leq 8r\sqrt{d}$.
\endproof

Proposition~\ref{prop:coordinate_wise_MOM} guarantees that starting from the coordinate-wise Median-of-Means we are off by  a $\sqrt{d}$ proportional factor from the optimal rate $r$. This will play a key role to analyze the number of steps we need to reach $\mu$ within the optimal rate $r$. Indeed, if we prove a geometric decay of the distance to $\mu$ along the descent step then only $\log d$ steps (up to a mutliplicative constants) would be enough to reach $\mu$ by a distance at most of the order of $r$.  \\

Let us now specify the step size we use at each iteration. At the current point $x_c$ we compute a top eigenvector $v_1$ of an approximating solution $M$ to \eqref{formule1} (i.e. $M$ such that $h_{x_c}(M)\geq (\beta\norm{x_c-\mu}_2 + 8r)^2$ for some $1/\sqrt{2}\leq \beta\leq 1$). Next iteration is $x_{c+1} = x_c - \theta_c v_1$ where the step size is
\begin{equation}\label{eq:step_size_def}
  \theta_c = - {\rm Med}\left(\inr{\bar X_k - x_c, v_1}:k\in[K]\right).
  \end{equation} In particular, since $\theta_c v_1$ does not depend  on the sign of $v_1$ (the product $\theta_c v_1$ is the same if we replace $v_1$ by $-v_1$), we do not care which top eigenvector of $M$ we choose.

  Let us now prove a geometric decay of the algorithm while $x_c$ is far from $\mu$. Again, this result is proved on the event $\cE$. 

  \begin{Proposition}\label{prop:geometric_decay}
  On the event $\cE$, the following holds. Let $x_c\in\bR^d$ (be the current point of the algorithm). Assume that $M$ is an approximating solution of \eqref{formule1}: $M$ is such that $h_{x_c}(M)\geq (\beta\norm{x_c-\mu}_2 + 8r)^2$ for some $0.78\leq \beta\leq 1$ and let $v_1$ be one of its top eigenvector. Then, we have 
\begin{equation*}
 \norm{x_{c+1} - \mu}_2^2\leq 0.8 \norm{x_c - \mu}_2^2 + 64r^2
\end{equation*} when $x_{c+1} = x_c - \theta_c v_1$ for $\theta_c$ defined in \eqref{eq:step_size_def}.  
  \end{Proposition}  

\beginproof Let us assume that the event $\cE$ holds up to the end of the proof. Let  $M$ be an approximating solution to \eqref{formule1}   such that $h_{x_c}(M)\geq (\beta\norm{x_c-\mu}_2 + 8r)^2$ for some $0.78\leq \beta\leq 1$ and let $v_1$ be a top eigenvector of $M$.

In direction $v_1$, there are at least $9K/10$ blocks such that $|\inr{\bar X_k - \mu, v_1}|\leq 8r$ hence on these blocks we also have
\begin{equation}\label{eq:prop_theta_c}
 |\theta_c - \inr{x_c-\mu,v_1}| = |{\rm Med}\left(\inr{\mu-\bar X_k,v_1}:k\in[K]\right)|\leq {\rm Med}\left(|\inr{\mu-\bar X_k,v_1}|:k\in[K]\right)\leq 8r.
 \end{equation}

Let $v= (\mu-x_c)/\norm{\mu-x_c}_2$ denote the optimal normalized descent direction.  We write $v = \lambda_1 v_1 + \lambda_2 v_1^\perp$ where $v_1^\perp$ is a normalized orthogonal vector to $v_1$.  We have $\lambda_1^2+\lambda_2^2=1$ and it follows from Proposition~\ref{Prop:direction} that $|\lambda_1| = |\inr{v_1, v}|>\sqrt{2\beta^2-1}$.  We conclude that
\begin{align*}
\norm{x_{c+1} - \mu}_2^2 & = \norm{x_c-\mu - \theta_c v_1}_2^2 = \norm{(\inr{x_c-\mu, v_1} -\theta_c)v_1  + \inr{x_c-\mu, v_1^\perp} v_1^\perp}_2^2\\ 
&= (\inr{x_c-\mu, v_1} -\theta_c)^2 +\inr{x_c-\mu, v_1^\perp}^2 \leq (8 r)^2 + \lambda_2^2 \norm{x_c-\mu}_2^2 
\end{align*} As $\lambda_2^2 = 1- \lambda_1^2<2-2\beta^2 < 0.8$ we get the result.
\endproof

We now have almost all the building blocks to fully characterize the algorithm. The last and final step is to find a stopping rule. The idea we use to design such a rule is based on Proposition~\ref{prop:geometric_decay}: we know that when the current point $x_c$ is not in a $\ell_2^d$-neighborhood of $\mu$ with a radius of the order of $r$ then the $\ell_2^d$-distance between the next iteration $x_{c+1}$ and $\mu$ should be less than $\sqrt{0.81}$ times the $\ell_2^d$-distance between $x_c$  and $\mu$. We therefore have a geometric decay of the distance to $\mu$ along the iterations until we reach $\mu$ in a $\ell_2^d$-neighborhood of radius proportional to $r$. Starting from the coordinate-wise median(-of-means) which is in a $8\sqrt{d}r$ neighborhood of $\mu$, we only have to do $\log(8\sqrt{d})/\log(1/\sqrt{0.81})$ iterations to output a current point which is $r$-close to $\mu$ w.r.t. the $\ell_2^d$-norm (see Proposition~\ref{prop:coordinate_wise_MOM}).

We are now in a position to write an ``almost final'' pseudo-code of our algorithm. In the next section, we will dive a bit deeper in this pseudo-code (and in particular on the covering SDP algorithm used to construct an approximating solution to \eqref{formule1}) in order to provide a final pseudo-code together with its total running time.

\vspace{0.7cm}
 \begin{algorithm}[H]\label{algo:almost_final}
\SetKwInOut{Input}{input}\SetKwInOut{Output}{output}\SetKw{Or}{or}
\SetKw{Return}{Return}
\Input{$X_1, \ldots, X_N$ and a number $K$ of blocks}
\Output{A robust subgaussian estimator of $\mu$}  
\BlankLine
Construct an equipartition $B_1\sqcup \cdots \sqcup B_K=\{1,\cdots,N\}$\\
Construct the $K$ empirical means $\bar{X}_k=(N/K)\sum_{i\in B_k}X_i, k\in[K]$\\
Compute $\hat\mu^{(0)}$ the coordinate-wise median-of-means and put $x_c \leftarrow \hat \mu^{(0)}$\\
\For{$T=1, 2, \cdots, \log(8\sqrt{d})/\log(1/\sqrt{0.81})$}{
Compute $M_c$ an approximating solution to \eqref{formule1} such that 
\begin{equation*}
h_{x_c}(M_c)\geq \left(0.78\norm{x_{c} - \mu}_2 + 8r\right)^2
\end{equation*}\\
Compute $v_1$ a top eigenvector of $M_{c}$\\
Compute a step size $\theta_{c} = -\Med\left(\inr{\bar X_k-x_{c},v_1}:k\in[K]\right)$\\
Update $x_c\leftarrow x_{c}-\theta_{c} v_1$\\}
\Return $x_{c}$
 \caption{``Almost final'' pseudo-code of the robust sub-gaussian estimator of $\mu$}
\end{algorithm}
\vspace{0.7cm}

Algorithm~\ref{algo:almost_final} is ``almost'' our final algorithm. There is one last step we need to check carefully: given a current point $x_c$ we need to find a way to construct $M_c$ satisfying ``$h_{x_c}(M_c)\geq \left(0.78\norm{x_{c} - \mu}_2 + 8r\right)^2$'' without knowing $r$ or $\mu$. This is the last issue we need to address in order to explain how step \textbf{5} from Algorithm~\ref{algo:almost_final} can be realized in a fully data-dependent way in a good time. This issue is answered in the next section together with the computation of its running time.

\section{Solving (approximatively) the SDP \eqref{formule1}}
\label{sec:SDP}
The aim of this section is to show that, on the event $\cE$, it is possible to construct in reasonnable time  a matrix $M_c$ such that  ``$h_{x_c}(M_c)\geq \left(0.78\norm{x_{c} - \mu}_2 + 8r\right)^2$'' without any extra information than the data. To that end we construct in an efficient way an approximation solution to the optimization problem \eqref{formule1} using covering SDP as in \cite{MR3909640}. The main result of this section is the following.

\begin{Theorem}\label{theo:approx_sol_sdp}Let $u\in\bN^*$.
On $ \cE$, for every $x_c\in\bR^d$ such that $\norm{x_c-\mu}_2\geq 800 r$, we can either compute, in time $\tilde\cO(K u d)$, with probability $> 1 -(1/10)^{u+5}/\sqrt{d}$ :
\begin{itemize}
\item A matrix $M_c$ such that 
\begin{equation*}
 h_{x_c}(M_c)\geq \left(0.78\norm{x_{c} - \mu}_2 + 8r\right)^2
\end{equation*}
\item Or directly a subgaussian estimate of $\mu$, using only the block means $\bar X_1, \ldots, \bar X_K$ as inputs. 
\end{itemize}
\end{Theorem}

Theorem~\ref{theo:approx_sol_sdp} answers the last issue raised at the end of Section~\ref{sec:proof_of_the_statistical_performance_in_theorem_theo:main} and provides the running time for step \textbf{5} of Algorithm~\ref{algo:almost_final}. It therefore concludes the statement that there exists a fully data-driven robust subgaussian algorithm for the estimation of a mean vector under the only Assumption~\ref{assum:first} (the total running time of Algorithm~\ref{algo:almost_final} is studied in Section~\ref{sec:final}).
\begin{Remark}
Theorem~\ref{theo:approx_sol_sdp} states that we either find an approximating solution $M_c$ to \eqref{formule1} or a good estimate of $\mu$ (at the current point $x_c$). As we will see in this section, this second case is degenerate as it is not the typical situation. 
\end{Remark}

We now turn to the proof of Theorem~\ref{theo:approx_sol_sdp}. It is decomposed into several lemmas adapted from  techniques developed by \cite{MR3909640} to approximately solve the semi-definite positive problem \eqref{formule1} in polynomial time. To that end, we first introduce the following covering SDP
 \begin{equation} \label{formule2} \tag{$C_\rho$}
\begin{aligned}
& \text{minimize}
& & \Tr(M^\prime)+ \norm{y^\prime}_1 \\
& \text{subject to}
& & M^\prime \succeq 0  , \ y^\prime \geq 0  , \\
&  &  & \forall k\in [K], \  \rho(\bar{X}_k- x_c)^\top M^\prime  (\bar{X}_k- x_c) + 9K/10 \  y_k^\prime \geq  1
\end{aligned}
\end{equation}where $\rho>0$ is some parameter that we will show how to fine-tune later. Then, we show that, for a good choice of $\rho$, we can turn a good approximation solution for \eqref{formule2} into a good approximation solution for \eqref{formule1}. 

We note $g(\rho)$ the optimal objective value of \eqref{formule2}. We begin with a first lemma that shows how to link the two optimization problems \eqref{formule1} and \eqref{formule2}. The proof can be found in Lemma 4.2 from \cite{MR3909640}. We adapt it here for our purpose.

\begin{Lemma}\label{SDP1}
Let $\rho>0$. From a feasible solution $(M^\prime, y^\prime)$ for \eqref{formule2} that achieves $\Tr(M^\prime)+ ||y^\prime||_1 \leq 1$, we can construct a feasible solution for \eqref{formule1} with objective value $\geq 1/\rho$ (and conversely).
\end{Lemma}
\beginproof
 We first note that the optimization problem \eqref{formule1} is equivalent to the following one:
 \begin{equation} \label{formule3} \tag{$\tilde E_{x_c}$}
\begin{aligned}
& \text{maximize}
& & z-\frac{\norm{y}_1}{9K/10} \\
& \text{subject to}
& & M \succeq 0 , \  \Tr(M)=1  , \ y \geq 0  , \ z\geq0 \\
&  &  & \forall k \in[K], \  (\bar{X}_k- x_c)^\top M (\bar{X}_k- x_c) +  \  y_k \geq  z
\end{aligned}
\end{equation}
Indeed, for a given $M\succeq0$ such that $\Tr(M)=1$, one can notice that the optimal value is achieved in \eqref{formule3} for $y_k=\max(0,z-(\bar{X}_k- x_c)^\top M (\bar{X}_k- x_c)), k\in[K]$ and $z=\mathcal{Q}_{9/10}\left( (\bar{X}_k- x_c)^\top M (\bar{X}_k- x_c) \right)$ the $9/10$-th quantile of $\{(\bar{X}_k- x_c)^\top M (\bar{X}_k- x_c):k\in[K]\}$, so that $z-\norm{y}_1/(9K/10)=h_{x_c}(M)$ which gives the equivalence between \eqref{formule1} and \eqref{formule3}. 

Then, once a feasible solution $(M^\prime, y^\prime)$ for \eqref{formule2} that achieves $\Tr(M^\prime)+ \norm{y^\prime}_1 \leq 1$ is obtained, by taking $M=M^\prime/\Tr(M^\prime)$, $z=1/(\rho \Tr(M^\prime))$ and $y=(9K/10)/(\rho \Tr(M^\prime)) y^\prime$, we get the desired result (and the converse follows from inverting those relations).
\endproof

From Lemma~\ref{SDP1}, it is enough to solve \eqref{formule2} -- for a good choice of $\rho$ --  to find a good approximating solution for \eqref{formule1}. It therefore remains to find such a good $\rho$. To do so, we rely on the next two lemmas. The first one is adapted from Lemma 4.3 in \cite{MR3909640}.  

\begin{Lemma}\label{SDP2}
For every $\rho>0$ and every $\alpha \in (0,1)$, $g( (1-\alpha)\rho) \geq g(\rho) \geq (1-\alpha) g( (1-\alpha)\rho)$.
\end{Lemma}

\beginproof A feasible pair $(M^\prime, y^\prime)$ for $(C_{(1-\alpha)\rho})$ is feasible for $(C_{\rho})$, which gives the first inequality. If   $(M^\prime, y^\prime)$ is a feasible pair for $(C_{\rho})$, then $(M^\prime/(1-\alpha), y^\prime/(1-\alpha))$ is a feasible pair for  $(C_{(1-\alpha)\rho})$, which gives the second inequality. 
\endproof

It follows from Lemma~\ref{SDP2}, that $g$ is continuous, non increasing, and (from Lemma \ref{SDP1}, using both sides of the implication, we have that $g(\rho)\leq1$ iff $1/\rho\geq OPT_{x_c}$) that $g(1/OPT_{x_c})=1$. So in order to find a good solution, we must find a $\rho$ such that $g(\rho)$ is as close to $1$ as possible. Unfortunately, we do not know how to solve \eqref{formule2} exactly for a given $\rho>0$, but we can compute efficiently a good approximation $(M^\prime, y^\prime)$ and a top eigenvector of $M^\prime$ thanks to the following result which can be found in \cite{PTZ12} and is detailed in  \cite{MR3909640} (see Section~4 and Remark~3.4). 

\begin{Lemma}\label{SDP4}[\cite{PTZ12}] Let $u\geq1$ be an integer. For every $\rho > 0$ and  every fixed $\eta > 0$, we can find with probability $>1-(1/10)^{u+10}/d$ a feasible solution to \eqref{formule2} that is $\eta$-close to the optimal, that is to say a feasible pair $(M^\prime, y^\prime)$ so that $\Tr(M^\prime)+\norm{y^\prime}_1 \leq (1+\eta)g(\rho)$ in time $ \tilde{\mathcal{O}}(uKd)$. Moreover, it is possible to find a top eigenvector of $M^\prime$ in $\tilde\cO(Kd)$.
\end{Lemma}

We compute $(u+ 3 \log(d)+10)$ times independently the (randomized) algorithm from \cite{PTZ12} that has a runtime of $ \tilde{\mathcal{O}}(Kd)$ and that outputs an     $\eta$-close feasible solution with probability $9/10$. By taking the largest of the output's objective value, we have an $\eta$-close feasible solution with probability $1-(1/10)^{u+3 \log(d)+10}$, in time $ \tilde{\mathcal{O}}(uKd)$, proving Lemma~\ref{SDP4}.  Let us call $\texttt{ALG}_\rho$ the algorithm from Lemma~\ref{SDP4}, that takes as input $((\bar{X}_k)_{k=1}^K, x_c, \rho, \eta, u)$ and returns a feasible pair $(M^\prime, y^\prime)$ for \eqref{formule2} satisfying $\Tr(M^\prime)+\norm{y^\prime}_1 \leq (1+\eta)g(\rho)$ in $\tilde{\mathcal{O}}(uKd)$, with probability $>1-(1/10)^{u+10}/d$. Next, in order to find a good $\rho$, we have to get some additional information on the function $g$. We will get it on the event $\cE$. 

\begin{Lemma}\label{SDP3}
On the event $\cal E$, for all $x_c\in\bR^d$, if $\norm{x_c-\mu}_2> 8r$ then
\begin{equation*}
g(\rho) \leq \frac{1}{\rho \ OPT_{x_c}} \left(1 +\rho  OPT_{x_c} \left(\frac{9(\norm{x_c-\mu}_2+8r)^2}{8(\norm{x_c-\mu}_2-8r)^2}-1\right)\right).
\end{equation*}

\end{Lemma}

\beginproof We use the same notation as in the  proof of Lemma~\ref{SDP1}. For any $\nu > 0$, we can choose a triplet $(z, y, M)$ feasible for \eqref{formule3} such that $z-\norm{y}_1/(9K/10) > OPT_{x_c}-\nu$. On the event $\cal E$, Lemma~\ref{DistanceLemme} yields $OPT_{x_c} > (8/9)(\norm{x_c-\mu}_2-8r)^2$ and we have from Corollary~\ref{coro:dual_isometry} that 
\begin{equation*}
z=\mathcal{Q}_{9/10}\left( (\bar{X}_k- x_c)^\top M (\bar{X}_k- x_c) \right) =\mathcal{Q}_{9/10}\left(\norm{M^{1/2}(\bar{X}_k- x_c)}_2\right)  \leq \left(\norm{M^{1/2}(x_c-\mu)}_2+8r\right)^2\leq (\norm{x_c-\mu}_2+8r)^2
\end{equation*}because $M\succeq0$ and $\Tr(M)=1$. Let $M^\prime= M/(\rho z),  y^\prime= y/[z(9K/10)]$. We have 
\begin{align*}
g(\rho) &\leq \Tr(M^\prime) + \norm{y^\prime}_1 \leq \frac{1+ \rho \norm{y}_1/(9K/10)}{\rho z}\\ &< \frac{1+\rho(z-OPT_{x_c} +\nu)}{\rho z } 
\leq \frac{1+\rho \nu +\rho OPT_{x_c}\left(\frac{ 9 (\norm{x_c-\mu}_2+8r)^2}{ 8(\norm{x_c-\mu}_2-8r)^2}-1\right) }{\rho (OPT_{x_c}-\nu)}.
\end{align*}
 By taking $\nu \rightarrow 0$, we get the result. 
\endproof

\textbf{Proof of Theorem~\ref{theo:approx_sol_sdp}.} Let us place ourselves on the event $\cE$ so that we can apply Lemma~\ref{SDP3}. Let $x_d\in\bR^d$ and assume that $\norm{x_c-\mu}_2 > 800 r$.   It follows from Lemma~\ref{SDP3} that $g(\rho) \leq 1/(\rho \ OPT_{x_c}) + 0.171$.  Therefore, if we can find a $\rho$ such that $g(\rho)\geq 1-\epsilon + 0.171$ for some  $0<\eps<1$, then necessarily $1/\rho \geq  OPT_{x_c} (1-\epsilon)$.  Let us take $\epsilon= 0.173$, and $\eta =0.0001$. Then if $\texttt{ALG}_\rho$ returns, a feasible pair $(M^\prime, y^\prime)$ for \eqref{formule2} so that $ 0.9981\leq \Tr(M^\prime)+\norm{y^\prime}_1 \leq 1$, then, since $0.9981> 1.0001\times 0.998 =(1+\eta)(1-\eps+0.171)$  we will know that, with probability $>1-(1/10)^{u+10}/d$, 
\begin{equation*}
(1+\eta)g(\rho)\geq \Tr(M^\prime)+\norm{y^\prime}_1\geq (1+\eta)(1-\eps+0.171)
\end{equation*}hence $1/\rho \geq OPT_{x_c}(1-\epsilon)$, and by Lemma \ref{SDP1}, we can  construct a feasible solution $M_c$ for \eqref{formule1} with objective value satisfying $h_{x_c}(M_c)\geq  OPT_{x_c}(1-\epsilon)$. Next, using Lemma~\ref{DistanceLemme}, we obtain that when $\norm{x_{c} - \mu}_2\geq 800r$
\begin{equation*}
h_{x_c}(M_c)\geq  OPT_{x_c}(1-\epsilon)\geq(1-\epsilon)(8/9)\left(\norm{x_{c} - \mu}_2 - 8r\right)^2\geq  \left(0.78\norm{x_{c} - \mu}_2 + 8r\right)^2
\end{equation*}for $\eps= 0.173$, solving step~\textbf{5} from Algorithm~\ref{algo:almost_final}. 

Therefore, it only remains to show how to find a $\rho$ such that $\texttt{ALG}_\rho$ returns a pair $(M^\prime, y^\prime)$ (feasible for \eqref{formule2}) satisfying $ 0.9981\leq \Tr(M^\prime)+\norm{y^\prime}_1 \leq 1$. We do it first by assuming that we have access to an initial $\rho_0$ such that $\texttt{ALG}_{\rho_0}$ returns a feasible pair $(M^\prime, y^\prime)$ for \eqref{formule2} (for $\rho=\rho_0$) so that  $\Tr(M^\prime) + \norm{y^\prime}_1\leq 1$ and to a maximal number $T$ of iterations  (we will also see later how to choose such $\rho_0$ and $T$). The following algorithm (which is a binary search) taking as input $(\bar{X}_1, \ldots, \bar{X}_K, x_c, \rho_0,u,  T)$ returns a feasible pair $(M^\prime, y^\prime)$ for \eqref{formule2} so that $ 0.9981\leq \Tr(M^\prime)+\norm{y^\prime}_1 \leq 1$ (when $T$ is large enough). This is simply due to the fact that $g$ is continuous, non increasing, $g(0)=10/9>1$ and $g(\rho)\leq2/8$ when $\rho\to+\infty$ and $\norm{x_c-\mu}_2>800r$ (because of Lemma~\ref{SDP3}). For this to work, we need that for each iteration, $\texttt{ALG}_{\rho}$ returns a feasible pair $(M^\prime, y^\prime)$ for \eqref{formule2} (for $\rho=\rho_0$) so that  $\Tr(M^\prime) + \norm{y^\prime}_1\leq (1+ 0.0001) g(\rho)$. We will suppose that it is the case for the rest of the proof. By union bound, this happens with probability at least $>1- T (1/10)^{u+10}/d$

\vspace{0.7cm}

\begin{algorithm}[H]\label{algo:binarySearch}
\SetKwInOut{Input}{input}\SetKwInOut{Output}{output}\SetKw{Or}{or}
\SetKw{Return}{Return}
\Input{$\bar{X}_1, \ldots, \bar{X}_K$, $x_c$, $\rho_0$,u,  $T$}
\Output{A feasible pair $(M^\prime, y^\prime)$ for \eqref{formule2} satisfying $ 0.9981\leq \Tr(M^\prime)+\norm{y^\prime}_1 \leq 1$}  
\BlankLine

$\rho_m \leftarrow 0$, $\rho_M \leftarrow \rho_0$, $V \leftarrow $ $\texttt{ALG}_{\rho_0}(u)$ , $i \leftarrow 0 $\\
\While{$ V  \notin[0,9981,1] $ and $i <T$}{\If{$V<0,9981$}{$\rho_M \leftarrow (\rho_M+\rho_m)/2$} \Else{$\rho_m \leftarrow (\rho_M+\rho_m)/2$}
$V \leftarrow objective(\texttt{ALG}_{\frac{\rho_m+\rho_M}{2}}(u))$ , $i \leftarrow i+1$}

\Return $\texttt{ALG}_{\frac{\rho_m+\rho_M}{2}(u)}$
 \caption{The \texttt{BinarySearch} algorithm to find a $\rho$ so that $\texttt{ALG}_\rho$ returns a pair $(M^\prime, y^\prime)$ (feasible for \eqref{formule2}) satisfying $ 0.9981\leq \Tr(M^\prime)+\norm{y^\prime}_1 \leq 1$.}
\end{algorithm}
\vspace{0.7cm}

If we can find a $\rho_0$ (such that $\texttt{ALG}_{\rho_0}$ returns a feasible pair $(M^\prime, y^\prime)$ for \eqref{formule2} so that  $\Tr(M^\prime) + \norm{y^\prime}_1\leq 1$) and a large enough number of iterations $T$ in \texttt{BinarySerach}, Algorithm \ref{algo:binarySearch} returns a feasible pair $(M^\prime, y^\prime)$ for \eqref{formule2} from which we can construct an approximating solution $M_c$ for \eqref{formule1} with objective value $h_{x_c}(M_c)$ larger than $\left(0.78\norm{x_{c} - \mu}_2 + 8r\right)^2$ whenever $\norm{x_c-\mu}_2 \geq 800 r$. This is exactly what we expect in step \textbf{5} of Algorithm~\ref{algo:almost_final}. Next, the last and final step that remains to be explained is to show how one can get such a $\rho_0$ and $T$ using only the block means $(\bar X_k)_{k=1}^K$ in $\tilde\cO(Nd+uKd)$.

Let us consider $\hat\mu^{(0)}$ the coordinate-wise median(-of-means) and let us define $\delta= \Med( \norm{\bar{X}_k-\hat\mu^{(0)}}_2:k\in[K])$ -- both  quantities can be computed in $\tilde{\mathcal{O}}(Kd)$. On the event $\cE$, it follows from Corollary~\ref{coro:dual_isometry}  (for  $M=I_d/d$) and Proposition~\ref{prop:coordinate_wise_MOM} that $\delta \leq 16 \sqrt{d} \times r$.  So if one takes $\rho_0= d/ \delta^2 \geq 1/[(16)^2r^2]$, and if $\norm{x_c-\mu}_2 > 800 r$, Lemma~\ref{DistanceLemme} and Lemma~\ref{SDP3} guarantee that $OPT_{x_c}\geq (8/9)\left(\norm{x_c-\mu}_2-8r\right)^2\geq (8/9)(792)^2r^2$ and so
\begin{equation*}
g(\rho_0)\leq \frac{1}{\rho\ OPT_{x_c}} + 0.171 \leq \frac{16^2}{(8/9)(792)^2} + 0.171<0.18
\end{equation*}
 so $\texttt{ALG}_{\rho_0}\leq (1+\eta)g(\rho)<1.0001\times 0.18<1$ (for the same choice of $\eta=0.0001$). 

 Now we tackle the question of the number $T$ of iterations, which is crucial for the runtime. We know from Lemma~\ref{SDP2} and  Lemma~\ref{SDP3}  that the interval $I$ of all $\rho$'s such that $ 0.9981 \leq  objective(\texttt{ALG}_{\rho}) \leq 1$ is at least of size $0.001/OPT_{x_c}$ when $\norm{x_c-\mu}_2 > 800 r$. Indeed, since $g(\rho)\leq objective(\texttt{ALG}_{\rho})\leq (1+\eta)g(\rho)$,  if $\rho$ is such that $0.9981\leq g(\rho)\leq1/(1+\eta)$ then $ 0.9981 \leq  objective(\texttt{ALG}_{\rho}) \leq 1$. Now, if we let $\rho_1>0$ and $0<\alpha<1$ be such that $g(\rho_1)=0.9981$ and $g((1-\alpha)\rho_1)=1/(1+\eta)$ the interval $I$ is at least of size $\alpha \rho_1$. Moreover, from Lemma~\ref{SDP2} we have $1/(1+\eta)\leq g((1-\alpha)\rho_1)\leq g(\rho_1)/(1-\alpha)$ and so $0.9981=g(\rho_1)\geq (1-\alpha)/(1+\eta)$, i.e. $\alpha\geq 1-0.9981(1+\eta)>0.001$. Finally, since $g(\rho_1)\leq 1$, $g(1/OPT_{x_c})=1$ and $g$ is non-increasing, we conclude that $\rho_1\geq 1/OPT_{x_c}$ and so the length of $I$ is at least $\alpha \rho_1\geq 0.001/OPT_{x_c}$.

  So, in the case where $\norm{x_c-\mu}_2>800r$,  $\log_2(\rho_0 \times OPT_{x_c}/0.001)$ iterations are enough to insure that \texttt{BinarySearch} outputs $(M^\prime, y^\prime)$ (from $\texttt{ALG}_{\rho}$ for a well-chosen $\rho$) feasible for \eqref{formule2} and such that $ 0.9981\leq \Tr(M^\prime)+\norm{y^\prime}_1 \leq 1$.  Moreover, on the event $\cE$ it is possible to show that for all iterations $x_c$ along the algorithm we have $\norm{x_c-\mu}_2 <  C \sqrt{d}r$ for a constant $C \leq 800$ (we may take that as an induction hypothesis for the firsts iterates $x_c$, and the proof of Theorem~\ref{theo:main} below in Section~\ref{sec:final} shows that it will still holds for $x_{c+1}$). So if $\delta > r/d$ then $\rho_0 <d^3/r^2 $, and since $OPT_{x_c} < (C^2 d +8) r^2$ (this follows from Lemma~\ref{DistanceLemme}), the binary search ends in time $T = \log_2(\tilde C d^4)$ with $\tilde C < 10^6$. \\

Thus, if the binary search has not ended in that time, we have either $\delta < r/d$ (which is a degenerate case) or $\norm{x_c-\mu}_2 < 800 r$ (or both). If $\norm{x_c-\mu}_2 > 800 r$ and $\delta < r/d$, then, taking $\rho_1= 1/  (d \delta)^2$, we have, by Lemma \ref{SDP3}, $ \texttt{ALG}_{\rho_1} < 1/2$. So, if we can not end our binary search in time $\log_2(\tilde C d^4)$, we compute $ \texttt{ALG}_{1/(d \delta)^2} $: if this gives something smaller than 1, that means that $1/(d \delta)^2 > 1/OPT_{x_c} \Rightarrow \delta < \sqrt{(C^2 d +8)} r/d < (C+1)r / \sqrt{d} $.  We notice that on $\cal E$, $\norm{\hat\mu^{(0)}-\mu}_2 < \delta+8r $, so if $ \texttt{ALG}_{1/(d \delta)^2} <1 $, then $\hat\mu^{(0)}$ is a good estimate for $\mu$. If on the contrary we have $ \texttt{ALG}_{\rho_1} > 1$, it means that $\norm{x_c-\mu}_2 < 800 r$, so we stop the algorithm and return $x_c$. \\

Let us write now in pseudo-code the procedure we just described. This is an algorithm, named \texttt{SolveSDP}, running in $\tilde\cO(Kud)$ which takes as inputs $\bar{X}_1, \ldots, \bar{X}_K$, $x_c$, $u$ and which outputs, on the event $\cE$, with probability $>1-\log(\tilde C d^4)(1/10)^{u+10}/d$, for every $x_c\in\bR^d$ such that $\norm{x_c-\mu}_2\geq 800 r$  either a matrix $M_c$ such that 
\begin{equation*}
 h_{x_c}(M_c)\geq \left(0.78\norm{x_{c} - \mu}_2 + 8r\right)^2
\end{equation*}
 or a subgaussian estimate of $\mu$. It therefore describes step~\textbf{5} from Algorithm~\ref{algo:almost_final}.

\vspace{0.3cm}
\begin{algorithm}[H]\label{algo:solveur}
\SetKwInOut{Input}{input}\SetKwInOut{Output}{output}\SetKw{Or}{or}
\SetKw{Return}{Return}
\Input{$\bar{X}_1, \ldots, \bar{X}_K$, $x_c$ and $u$}
\Output{A feasible solution for \eqref{formule1}}  
\BlankLine
Compute $\hat\mu^{(0)}$, compute $\delta$\\
$T \leftarrow  \log(\tilde C d^4)$, $\rho_0 \leftarrow d/ \delta^2$\\
$(M^\prime, y^\prime) \leftarrow$ BinarySearch($T$, $\rho_0$)\\
\If{$\Tr(M^\prime)+||y||_1 \in[0,9981,1]$ }{$M\leftarrow M^\prime/\Tr(M^\prime)$\\ \Return (True, $M$)}
\Else{\If{ $\texttt{ALG}_{1/(d \delta)^2} <1$}{\Return(False, $\hat\mu^{(0)}$) } \Else{\Return(False, $x_c$)}}

 \caption{SolveSDP}
\end{algorithm}

\begin{Remark}\label{rem:effects_blocks}[Two advantages of block means]
During the whole algorithm, we solve the program \eqref{formule2} up to a factor $(1+ \eta)$ where $ \eta $ is \emph{fixed} (here we take it equal to  $0.0001$). This differs crucially from the work of \cite{MR3909640} where $\eta$ depends on the fraction of outliers, which decreases the performance of the algorithm in Lemma \ref{SDP4}, the true runnnig time being $\tilde \cO(Kd/\text{Poly}(\eta))$. This is another advantages of using the mean blocks instead of the data themselves. Indeed, using blocks of data, we work with a constant fraction of corrupted blocks (we took it equal to $1/10$), therefore the approximation parameter used to approximately solved \eqref{formule2} can be taken equal to a constant (we took it equal to $\eta=0.0001$) unlike \cite{MR3909640} where $\eta$ depends on $\eps=|\cO|/N$. Taking the block means has therefore two advantages: a stochastic one, which is to exhibit a subgaussian behavior for $9K/10$ blocks even under a $L_2$-moment assumption and a computational one, which is to make the proportion of corrupted blocks constant.
\end{Remark}

 \section{The final algorithm and its computational cost: proof of Theorem~\ref{theo:main}.}
 \label{sec:final}
 We are now in a position to fully describe our robust subgaussian descent algorithm running in $\tilde \cO(Nd+uKd)$. One may check that  its construction is fully data-dependent, in particular, we do not need to know the value of $r$ or the proportion of outliers. 

 \vspace{0.7cm}
 \begin{algorithm}[H]\label{algo:final}
\SetKwInOut{Input}{input}\SetKwInOut{Output}{output}\SetKw{Or}{or}
\SetKw{Return}{Return}
\Input{$X_1, \ldots, X_N$ and $K\in[N]$ and $u\in\bN^*$}
\Output{A robust subgaussian estimator of $\mu$}
\BlankLine
Construct an equipartition $B_1\sqcup \cdots \sqcup B_K=\{1,\cdots,N\}$\\
Construct the $K$ empirical means $\bar{X}_k=(N/K)\sum_{i\in B_k}X_i, k\in[K]$\\
Compute $\hat\mu^{(0)}$ the coordinate-wise median\\
$x_c \leftarrow \hat \mu^{(0)}$, Bool $\leftarrow$ True, $T \leftarrow 0$ \\
\While{Bool and $T <\log(8\sqrt{d})/\log(1/0.81)$}{Bool, $A$  $\leftarrow$SolveSDP($\bar{X}_1, \ldots, \bar{X}_K$, $x_c$)\\
\If{Bool}{$M \leftarrow A$ \\ Compute $v_1$ a top eigenvector of $M_{c}$\\
Compute a step size $\theta_{c} = -\Med\left(\inr{\bar X_k-x_{c},v_1}:k\in[K]\right)$\\
Update $x_{c}\leftarrow x_{c}-\theta_{c} v_1$ \\ $T\leftarrow T+1$\\ }\Else{$x_c \leftarrow A$}}
\Return $x_c$
 \caption{Final Algorithm: covSDPofMeans}
\end{algorithm}
\vspace{0.7cm}


\textbf{Proof of Theorem~\ref{theo:main}.} From Theorem~\ref{theo:approx_sol_sdp}, we know that on $\cE $, when, $\norm{x_c-\mu}_2> 800r$, we get, with probability $>1-(1/10)^{u+5}/\sqrt{d}$, an $M_c$ so that $ h_{x_c}(M_c)\geq \left(0.8\norm{x_{c} - \mu}_2 + 8r\right)^2$  (or directly a subgaussian estimate, in which case our work is done). Proposition \ref{prop:geometric_decay}, states that in that case  $\norm{x_{c+1} - \mu}_2^2\leq 0.8 \norm{x_c - \mu}_2^2 + 64r^2 \leq 0.81 \norm{x_c - \mu}_2^2$. So we have a geometric decays and Proposition \ref{prop:coordinate_wise_MOM} guarantees that our starting point is at most $8\sqrt{d} r$ far away from the mean so that in at most $\log(8\sqrt{d})/\log(1/0.81))$ steps the algorithm outputs its current point which is $r$-close to $\mu$, with probability $>1-(1/10)^{u+5} \log(8\sqrt{d})/(\log(1/0.81))\sqrt{d})>1-(1/10)^u$ (by union bound).

The last thing to do is to control what happens when $\norm{x_c-\mu}_2< 800r$. Then, we have no guarantees on $v_1$, but using the similar argument as in the proof of Proposition~\ref{prop:geometric_decay} we know that 
\begin{equation}\label{eq:prop_theta_c}
 |\theta_c - \inr{x_c-\mu,v_1}| = |{\rm Med}\left(\inr{\mu-\bar X_k,v_1}:k\in[K]\right)|\leq {\rm Med}\left(|\inr{\mu-\bar X_k,v_1}|:k\in[K]\right)\leq 8r
 \end{equation} and (for some $v_1^\perp$ a normalized orthogonal vector to $v_1$)  
\begin{align*}
\norm{x_{c+1} - \mu}_2^2 & = \norm{x_c-\mu - \theta_c v_1}_2^2 = \norm{(\inr{x_c-\mu, v_1} -\theta_c)v_1  + \inr{x_c-\mu, v_1^\perp} v_1^\perp}_2^2\\ 
&= (\inr{x_c-\mu, v_1} -\theta_c)^2 +\inr{x_c-\mu, v_1^\perp}^2 \leq (8 r)^2 +  \norm{x_c-\mu}_2^2 .
\end{align*} Hence,  $\norm{x_{c+1} - \mu}_2 \leq (8 r) +  \norm{x_c-\mu}_2 $. Therefore, in the worst case scenario where $\norm{x_c-\mu}_2 >800r $ at the last iteration, the algorithm outputs the next iteration $\hat \mu_K = x_{c+1}$ so that $\norm{\hat \mu_K - \mu}_2\leq 808r$.


We end this proof with the computation of the running time of Algorithm~\ref{algo:final}. We detail the computation cost for each line of Algorithm~\ref{algo:final}: line~\textbf{1} cost $N$, line~\textbf{2} costs $Nd$, line~\textbf{3} costs $\cO(d K \log(K))$. The while loop in line~\textbf{5} is running at least $\log d$ times (up to constant) so that the computational cost of all remaining lines of Algorithm~\ref{algo:final} are at worst to be multiplied by $\log d$. Line~\textbf{6} costs $\log(\tilde C d^4)$ steps, each of cost $\tilde \cO(Kud)$ (that comes from Lemma \ref{SDP4}). Line~\textbf{9} can be computed in $\tilde\cO(Nd)$ thanks to Lemma~\ref{SDP4}. Finally,  line~\textbf{10} costs $\cO(Kd)$. Other lines take time at most $d$. We thus recover the running time announced in Theorem~\ref{theo:main}. 
\endproof

\section{Adaptive choice of $K$} 
\label{sec:adaptive_choice_of_}
Given a number of blocks $K\in\{1, \ldots, N\}$, a parameter $u\geq1$ (so that  the covering SDPs from \cite{PTZ12} (used in Lemma~\ref{SDP4}) is ran $u+3 \log d+10$ times) and the dataset $\{X_1, \ldots, X_N\}$, Algorithm~\ref{algo:final} returns a vector $\hat \mu_K$ in $\bR^d$ and Theorem~\ref{theo:main} insures that $\hat \mu_K$ estimates the true mean $\mu$ at the subgaussian rate \eqref{eq:intro_subgaus_rate} with large probability as long as $K\geq 300|\cO|$. As a consequence, we have certified statistical guarantees for $\mu_K$ only when some a priori knowledge on the number $|\cO|$ of outliers is provided (such as ``the corruption of this database is less than $5\%$'' ) or if we choose $K$ like $N$- but, in this later case the rate  \eqref{eq:intro_subgaus_rate} may be too pessimistic. The aim of this section is to overcome this issue by constructing a procedure which can automatically adapt to the number of outliers. The resulting procedure satisfies the same statistical bounds as $\mu_K$ for all $K\geq 300|\cO|$ without knowing $|\cO|$ (up to constants). 

The adaptation method we use is based on the Lepski method \cite{MR1091202,MR1147167} which is another tool used by the ``MOM community'' since \cite{lugosi2019sub}. The price we pay for this adaptation is the a priori knowledge of the rate \eqref{eq:intro_subgaus_rate} for all $K$ which means that we know in advance $\Tr(\Sigma)$ and $\norm{\Sigma}_{op}$ -- this is for instance the case when it is known that $\Sigma$ is the identity matrix $I_d$. Of course, one can design robust estimators for $\Tr(\Sigma)$ (see \cite{Jules_Guillaume_1}) and $\norm{\Sigma}_{op}$ but this requires stronger assumptions that we want to avoid at this stage.  

Lepski's method proceeds as follows. We set for all $K\in\{1, \ldots, N\}$ and all $j\in\{0,1,\ldots, \log_2N\}$
\begin{equation*}
r_K^* = 808\left(1200\sqrt{\frac{\Tr(\Sigma)}{N}} + \sqrt{\frac{1200\norm{\Sigma}_{op}K}{N}}\right) \mbox{ and } r^{(j)} = r^*_{\lceil N/2^j\rceil}
\end{equation*}the rate of convergence from Theorem~\ref{theo:main}. For a given parameter $u_j\in\bN^*$, we construct from Algorithm~\ref{algo:final}
\begin{equation}\label{eq:def_esti_lepski}
\hat \mu^{(j)}\leftarrow covSDPofMeans(X_1,\ldots,X_N, K=\lceil N/2^j\rceil, u=u_j).
\end{equation}
Classical Lepski's method considers the largest $J$ such that $\cap_{j=0}^J B_2(\hat \mu^{(j)}, r^{(j)})$ is none empty and then take any point $\hat \mu$ in this none empty intersection. Standard analysis of Lepski's method shows that $\hat \mu$ estimates $\mu$ at the rate $r_K^*$ (up to an absolute constant) simultaneously for all $K\in \{300 |\cO|, \ldots, N\}$ without knowing $|\cO|$. Given that checking that the intersection of several $\ell_2^d$-balls may not be straigtforward, we use a slightly modified version of Lepski's method as described in the following algorithm.


 \vspace{0.7cm}
 \begin{algorithm}[H]\label{algo:lepski}
\SetKwInOut{Input}{input}\SetKwInOut{Output}{output}\SetKwInOut{Init}{init}\SetKw{Or}{or}
\SetKw{Return}{Return}
\Input{$X_1, \ldots, X_N$ and $\{u_j:j=0,1,2,\ldots,\log_2 N\}\subset \bN^*$}
\Output{A robust subgaussian estimator of $\mu$ with adaptive choice of $K$}
\Init{$J=0$ and $\hat \mu^{(0)} = covSDPofMeans(X_1,\ldots,X_N, K=N, u=u_0)$}
\BlankLine
\While{$\norm{\hat \mu^{(J)} - \hat \mu^{(j)}}_2\leq r^{(J)} + r^{(j)} , j= J-1,J-2,\ldots, 0$}{
$J\leftarrow J+1$\\
$\hat \mu^{(J)}\leftarrow covSDPofMeans(X_1,\ldots,X_N, K=\lceil N/2^J\rceil, u=u_J)$	
}
\Return $\hat \mu^{(J)}$
 \caption{Adaptive choice of $K$ in covSDPofMeans}
\end{algorithm}
\vspace{0.7cm}

Unlike for the traditional Lepski's method we check that $\hat\mu^{(J)}$ is in $\cap_{j=0}^{J-1} B_2(\hat \mu^{(j)}, r^{(J)} +r^{(j)})$ instead of checking that $\cap_{j=0}^J B_2(\hat \mu^{(j)}, r^{(j)})$ is none empty -- this simplifies the adaptation step. It is also possible to speed up the whole procedure by constructing iteratively the block means. Indeed, given that we consider a dyadic grid for $K$, i.e. $K\in\{N,\lceil N/2\rceil,\lceil N/4\rceil, \ldots\}$, for all $j\in\bN$, we can construct the block means $\{\bar X_k^{(j+1)},k=1,\ldots, \lceil N/2^{j+1}\rceil\}$ at step $K=\lceil N/2^{j+1}\rceil$ using the block means from the previous step $K=\lceil N/2^{j}\rceil$  by simply averaging two successive block means: $\bar X_k^{(j+1)}\leftarrow (\bar X_{2k}^{(j)} + \bar X_{2k+1}^{(j)})/2$.

Let us now turn to the statistical analysis of the output $\hat \mu^{(\hat J)}$ from Algorithm~\ref{algo:lepski} where 
\begin{equation*}
 \hat J =\max\left(J\in\{0,1,\ldots, \log_2 N\}: \hat \mu^{(J)}\in\cap_{j=0}^{J-1} B_2(\hat \mu^{(j)}, r^{(J)} + r^{(j)})\right).
 \end{equation*} 

\begin{Theorem}\label{theo:lepski}
Let $\{u_j:j=0,1,2,\ldots,\log_2 N\}\subset \bN^*$ be the family of parameters used to construct the family of estimators  $\{\hat \mu^{(j)}, j=0, 1, \ldots\}$ in Algorithm~\ref{algo:lepski} (see also \eqref{eq:def_esti_lepski}). For all $K\in\{600|\cO|, \ldots, N\}$, with probability at least 
\begin{equation}\label{eq:proba_lepski}
1-2\exp(-K/360000)-\sum_{j=0}^{\log_2(N/(K-1))}(1/10)^{u_j}
\end{equation} the output $\hat \mu^{(\hat J)}$ of Algorithm~\ref{algo:lepski} is such that $\norm{\hat \mu^{(\hat J)} - \mu}_2\leq 3 r^*_K$.
\end{Theorem}

\beginproof
 For all $j\in\{0,1,\ldots, \log_2 N\}$ denote by $\cE_j$ the event onto which Theorem~\ref{theo:main} is valid for $K=\lceil N/2^{j}\rceil$ and for $u=u_j$: that is on $\cE_j$, if $\lceil N/2^{j}\rceil\geq 300|\cO|$, $\norm{\hat \mu^{(j)}-\mu}_2\leq r^{(j)}$ and $\bP[\cE_j]\geq 1-\exp(-\lceil N/2^{j}\rceil/180000)-(1/10)^{u_j}$.  Let $K\in\{600|\cO|, \ldots, N\}$ and $J\in\{0,1,\ldots, \log_2 N\}$ be such that $\lceil N/2^{J}\rceil\leq K <\lceil N/2^{J-1}\rceil$. On the event $\cap_{j=0}^J \cE_{j}$, we have $\norm{\hat \mu^{(j)}-\mu}_2\leq r^{(j)}$ for all $j=0,1,\ldots,J$, in particular, for all $j=0,1,\ldots, J-1$, $\norm{\hat \mu^{(J)} - \hat \mu^{(j)}}_2\leq r^{(J)} + r^{(j)}$ and so $\hat \mu^{(J)}\in\cap_{j=0}^{J-1} B_2(\hat \mu^{(j)}, r^{(J)} +r^{(j)})$. As a consequence $\hat J\geq J$ therefore $\norm{\hat \mu^{(\hat J)} - \hat \mu^{(J)}}_2\leq r^{(\hat J)} + r^{(J)}\leq 2 r^{(J)}\leq 2 r_K^*$. Finally, we have
 \begin{align*}
 \bP[\cap_{j=0}^J \cE_{j}] \geq 1-\sum_{j=0}^J \exp(-\lceil N/2^{j}\rceil/180000)-(1/10)^{u_j}\geq 1-2\exp(-K/360000)-\sum_{j=0}^{\log_2(N/(K-1))}(1/10)^{u_j}.
 \end{align*}
 \endproof

 We can see in Algorithm~\ref{algo:lepski}  that $\hat \mu^{(\hat J)}$ does not use any information on the number of outliers $|\cO|$ for its construction but it can still estimate $\mu$ at the optimal rate $r^*_K$ for all deviation parameters $K$ in $\{600|\cO|, \ldots, N\}$. The maximum total running time of Algorithm~\ref{algo:lepski} is achieved when $\hat J = \log_2N$; in that case, it is at most $\tilde\cO(Nd + \sum_{j=0}^{\log_2N} \lceil N/2^j\rceil u_j d)$. In particular, if one chooses $u_j = 2^j$ for all $j=0,1,\ldots, \log_2N$ then the total running time for the construction of $\hat \mu^{(\hat J)}$ is nearly-linear $\tilde \cO(Nd)$. For this choice of $u_j$, the probability deviation in \eqref{eq:proba_lepski} is constant and so one should choose the smallest possible $K$ allowed in Theorem~\ref{theo:lepski}, that is $K=600|\cO|$. Let us write formally this result.

 \begin{Corollary}\label{coro:lepski}
 If one takes $u_j=2^j$ for all $j=0,1,\ldots, \log_2N$ in Algorithm~\ref{algo:lepski} then, in nearly-linear time $\tilde\cO(Nd)$, with probability at least $1-2\exp(-600|\cO|/360000)-1/11$, the output $\hat \mu^{(\hat J)}$ from Algorithm~\ref{algo:lepski} satisfies
 \begin{equation}\label{eq:rate_coro_lepski}
  \norm{\hat\mu^{(\hat J)} - \mu}_2\leq 2r^*_{600|\cO|}=1616\left(1200\sqrt{\frac{\Tr(\Sigma)}{N}} + 850\sqrt{\frac{\norm{\Sigma}_{op}|\cO|}{N}}\right).
  \end{equation} 
 \end{Corollary}
In particular, considering the setup from Theorem~\ref{theo:diakonikolas}, if $|\cO| = \eps N$ for some $\eps\leq 1/600$ then the rate achieved by $\hat\mu^{(\hat J)}$ in Corollary~\ref{coro:lepski} is of the order of 
\begin{equation*}
\sqrt{\frac{\Tr(\Sigma)}{N}} + \sqrt{\norm{\Sigma}_{op}\eps}
\end{equation*}which is like $\sqrt{\norm{\Sigma}_{op}\eps}$ when $N\geq (\Tr(\Sigma)/\norm{\Sigma}_{op})/\eps$. As a consequence, the result from Corollary~\ref{coro:lepski} improves the one from Theorem~\ref{theo:diakonikolas} by removing an extra $\log d$ factor in the sample complexity in the case considered in Theorem~\ref{theo:diakonikolas} that is when $\Sigma\preceq \sigma^2 I_d$.  Moreover, Corollary~\ref{coro:lepski} also shows that the sample complexity depends on the \textit{effective rank} $\Tr(\Sigma)/\norm{\Sigma}_{op}$ of $\Sigma$. This ratio can be much smaller than $d$ if the spectrum of $\Sigma$ decays sufficiently fast. Finally, Corollary~\ref{coro:lepski} also covers the case where the sample size $N$ is less than the sample complexity -- that is when $N\leq (\Tr(\Sigma)/\norm{\Sigma}_{op})/\eps$. In that case, the estimation rate is given by $\sqrt{\Tr(\Sigma)/N}$ which is the complexity coming from the estimation of $\mu$ in the none corrupted case. As a consequence, Corollary~\ref{coro:lepski} exhibits a phase transition happening at $N\sim (\Tr(\Sigma)/\norm{\Sigma}_{op})/\eps$ above which corruption is the main source of estimation mistakes and below which corruption does not play any role.

Corollary~\ref{coro:lepski} covers the case where $\hat\mu^{(\hat J)}$ is computed in nearly-linear time and with statistical guarantees happening with constant probability. In the following final result, we show that $\hat\mu^{(\hat J)}$ can estimate $\mu$ at the optimal rate $r^*_K$ for all $K\geq 600|\cO|$ with a subgaussian deviation $1-2\exp(-K/360000)$ if we perform more iterations $u_j$ of the covering SDP from Lemma~\ref{SDP4}. The price we pay for this subgaussian behavior of $\hat\mu^{(\hat J)}$  is on the total running time which goes from nearly-linear time $\tilde\cO(Nd)$ to $\tilde\cO(N^2d)$ by taking $u_j = \lceil N/2^j\rceil$ for $j=0,1,\ldots, \log_2N$ ($u_j=N$ would do as well). We write formally this statement in the next corollary which follows directly from Theorem~\ref{theo:lepski}.

\begin{Corollary}\label{coro:lepski_2_subgaussian}
 If one takes $u_j=\lceil N/2^j\rceil$ for all $j=0,1,\ldots, \log_2N$ in Algorithm~\ref{algo:lepski} then, in time $\tilde\cO(N^2d)$, for all $K\geq 600|\cO|$, with probability at least $1-4\exp(-K/360000)$, the output $\hat \mu^{(\hat J)}$ from Algorithm~\ref{algo:lepski} satisfies
 \begin{equation}\label{eq:rate_coro_lepski}
  \norm{\hat\mu^{(\hat J)} - \mu}_2\leq 2r^*_{K}=1616\left(1200\sqrt{\frac{\Tr(\Sigma)}{N}} + \sqrt{\frac{1200\norm{\Sigma}_{op}K}{N}}\right).
  \end{equation} 
\end{Corollary}
As a consequence $\hat \mu^{(\hat J)}$ is a subgaussian estimator of $\mu$ for all range of $K$ from $600|\cO|$ to $N$ which can handle up to $|\cO|$ outliers in the database (even when $|\cO|\sim N$) and that can be constructed in time $\tilde\cO(N^2d)$. It does not require any knowledge on $|\cO|$ for its construction.

\vspace{0.7cm}
\textbf{Acknowlegements:} We would like to thank Yeshwanth Cherapanamjeri, Ilias Diakonikolas, Yihe Dong, Nicolas Flammarion, Sam Hopkins and Jerry Li for helpful comments on our work.





\begin{footnotesize}
\bibliographystyle{plain}
\bibliography{biblio}

\begin{thebibliography}{10}

\bibitem{MR1688610}
Noga Alon, Yossi Matias, and Mario Szegedy.
\newblock The space complexity of approximating the frequency moments.
\newblock {\em J. Comput. System Sci.}, 58(1, part 2):137--147, 1999.
\newblock Twenty-eighth Annual ACM Symposium on the Theory of Computing
  (Philadelphia, PA, 1996).

\bibitem{MR3185193}
St\'{e}phane Boucheron, G\'{a}bor Lugosi, and Pascal Massart.
\newblock {\em Concentration inequalities}.
\newblock Oxford University Press, Oxford, 2013.
\newblock A nonasymptotic theory of independence, With a foreword by Michel
  Ledoux.

\bibitem{MR3124669}
S\'{e}bastien Bubeck, Nicol\`o Cesa-Bianchi, and G\'{a}bor Lugosi.
\newblock Bandits with heavy tail.
\newblock {\em IEEE Trans. Inform. Theory}, 59(11):7711--7717, 2013.

\bibitem{AIHPB_2012__48_4_1148_0}
Olivier Catoni.
\newblock Challenging the empirical mean and empirical variance: A deviation
  study.
\newblock {\em Annales de l'I.H.P. Probabilit\'es et statistiques},
  48(4):1148--1185, 2012.

\bibitem{CG}
Olivier Catoni and Ilaria Giulini.
\newblock Dimension-free pac-bayesian bounds for matrices, vectors, and linear
  least squares regression.
\newblock Technical report, CNRS and LSPM, 2017.

\bibitem{MR3845006}
Mengjie Chen, Chao Gao, and Zhao Ren.
\newblock Robust covariance and scatter matrix estimation under {H}uber's
  contamination model.
\newblock {\em Ann. Statist.}, 46(5):1932--1960, 2018.

\bibitem{MR3909640}
Yu~Cheng, Ilias Diakonikolas, and Rong Ge.
\newblock High-dimensional robust mean estimation in nearly-linear time.
\newblock In {\em Proceedings of the {T}hirtieth {A}nnual {ACM}-{SIAM}
  {S}ymposium on {D}iscrete {A}lgorithms}, pages 2755--2771. SIAM,
  Philadelphia, PA, 2019.

\bibitem{Bartlett19}
Yeshwanth Cherapanamjeri, Nicolas Flammarion, and Peter~L. Bartlett.
\newblock Fast mean estimation with sub-gaussian rates, 2019.

\bibitem{MR1666908}
V\'{\i}ctor~H. de~la Pe\~{n}a and Evarist Gin\'{e}.
\newblock {\em Decoupling}.
\newblock Probability and its Applications (New York). Springer-Verlag, New
  York, 1999.
\newblock From dependence to independence, Randomly stopped processes.
  $U$-statistics and processes. Martingales and beyond.

\bibitem{MR3576558}
Luc Devroye, Matthieu Lerasle, Gabor Lugosi, and Roberto~I. Oliveira.
\newblock Sub-{G}aussian mean estimators.
\newblock {\em Ann. Statist.}, 44(6):2695--2725, 2016.

\bibitem{MR3945261}
Ilias Diakonikolas, Gautam Kamath, Daniel Kane, Jerry Li, Ankur Moitra, and
  Alistair Stewart.
\newblock Robust {E}stimators in {H}igh-{D}imensions {W}ithout the
  {C}omputational {I}ntractability.
\newblock {\em SIAM J. Comput.}, 48(2):742--864, 2019.

\bibitem{MR3631028}
Ilias Diakonikolas, Gautam Kamath, Daniel~M. Kane, Jerry Li, Ankur Moitra, and
  Alistair Stewart.
\newblock Robust estimators in high dimensions without the computational
  intractability.
\newblock In {\em 57th {A}nnual {IEEE} {S}ymposium on {F}oundations of
  {C}omputer {S}cience---{FOCS} 2016}, pages 655--664. IEEE Computer Soc., Los
  Alamitos, CA, 2016.

\bibitem{diakonikolas2016robust}
Ilias Diakonikolas, Gautam Kamath, Daniel~M Kane, Jerry Li, Ankur Moitra, and
  Alistair Stewart.
\newblock Robust estimators in high dimensions without the computational
  intractability.
\newblock In {\em Foundations of Computer Science (FOCS), 2016 IEEE 57th Annual
  Symposium on}, pages 655--664. IEEE, 2016.

\bibitem{diakonikolas2018robustly}
Ilias Diakonikolas, Gautam Kamath, Daniel~M Kane, Jerry Li, Ankur Moitra, and
  Alistair Stewart.
\newblock Robustly learning a gaussian: Getting optimal error, efficiently.
\newblock In {\em Proceedings of the Twenty-Ninth Annual ACM-SIAM Symposium on
  Discrete Algorithms}, pages 2683--2702. Society for Industrial and Applied
  Mathematics, 2018.

\bibitem{MR3826316}
Ilias Diakonikolas, Daniel~M. Kane, and Alistair Stewart.
\newblock List-decodable robust mean estimation and learning mixtures of
  spherical {G}aussians.
\newblock In {\em S{TOC}'18---{P}roceedings of the 50th {A}nnual {ACM} {SIGACT}
  {S}ymposium on {T}heory of {C}omputing}, pages 1047--1060. ACM, New York,
  2018.

\bibitem{MR3909639}
Ilias Diakonikolas, Weihao Kong, and Alistair Stewart.
\newblock Efficient algorithms and lower bounds for robust linear regression.
\newblock In {\em Proceedings of the {T}hirtieth {A}nnual {ACM}-{SIAM}
  {S}ymposium on {D}iscrete {A}lgorithms}, pages 2745--2754. SIAM,
  Philadelphia, PA, 2019.

\bibitem{MR1193313}
David~L. Donoho and Miriam Gasko.
\newblock Breakdown properties of location estimates based on halfspace depth
  and projected outlyingness.
\newblock {\em Ann. Statist.}, 20(4):1803--1827, 1992.

\bibitem{MR0301858}
Frank~R. Hampel.
\newblock A general qualitative definition of robustness.
\newblock {\em Ann. Math. Statist.}, 42:1887--1896, 1971.

\bibitem{MR0359096}
Frank~R. Hampel.
\newblock Robust estimation: a condensed partial survey.
\newblock {\em Z. Wahrscheinlichkeitstheorie und Verw. Gebiete}, 27:87--104,
  1973.

\bibitem{hopkins2018sub}
Samuel~B Hopkins.
\newblock Sub-gaussian mean estimation in polynomial time.
\newblock {\em arXiv preprint arXiv:1809.07425}, 2018.

\bibitem{MR0161415}
Peter~J. Huber.
\newblock Robust estimation of a location parameter.
\newblock {\em Ann. Math. Statist.}, 35:73--101, 1964.

\bibitem{MR2488795}
Peter~J. Huber and Elvezio~M. Ronchetti.
\newblock {\em Robust statistics}.
\newblock Wiley Series in Probability and Statistics. John Wiley \& Sons, Inc.,
  Hoboken, NJ, second edition, 2009.

\bibitem{MR855970}
Mark~R. Jerrum, Leslie~G. Valiant, and Vijay~V. Vazirani.
\newblock Random generation of combinatorial structures from a uniform
  distribution.
\newblock {\em Theoret. Comput. Sci.}, 43(2-3):169--188, 1986.

\bibitem{Led01}
Michel Ledoux.
\newblock {\em The concentration of measure phenomenon}, volume~89 of {\em
  Mathematical Surveys and Monographs}.
\newblock American Mathematical Society, Providence, RI, 2001.

\bibitem{MR2814399}
Michel Ledoux and Michel Talagrand.
\newblock {\em Probability in {B}anach spaces}.
\newblock Classics in Mathematics. Springer-Verlag, Berlin, 2011.
\newblock Isoperimetry and processes, Reprint of the 1991 edition.

\bibitem{LO}
M.~Lerasle and R.~Oliveira.
\newblock Robust empirical mean estimators.
\newblock Technical report, IMPA and CNRS, 2011.

\bibitem{lugosi2019sub}
G{\'a}bor Lugosi, Shahar Mendelson, et~al.
\newblock Sub-gaussian estimators of the mean of a random vector.
\newblock {\em The Annals of Statistics}, 47(2):783--794, 2019.

\bibitem{LMSL}
Z.~Szabo M.~Lerasle, T.~Matthieu and G.~Lecué.
\newblock Monk – outliers-robust mean embedding estimation by
  median-of-means.
\newblock Technical report, CNRS, University of Paris 11, Ecole Polytechnique
  and CREST, 2017.

\bibitem{MS}
S~Minsker and N.~Strawn.
\newblock Distributed statistical estimation and rates of convergence in normal
  approximation.
\newblock Technical report, arXiv: 1704.02658, 2017.

\bibitem{minsker2015geometric}
Stanislav Minsker.
\newblock Geometric median and robust estimation in banach spaces.
\newblock {\em Bernoulli}, 21(4):2308--2335, 2015.

\bibitem{MR3851758}
Stanislav Minsker.
\newblock Sub-{G}aussian estimators of the mean of a random matrix with
  heavy-tailed entries.
\newblock {\em Ann. Statist.}, 46(6A):2871--2903, 2018.

\bibitem{MR702836}
A.~S. Nemirovsky and D.~B.~and Yudin.
\newblock {\em Problem complexity and method efficiency in optimization}.
\newblock A Wiley-Interscience Publication. John Wiley \& Sons, Inc., New York,
  1983.
\newblock Translated from the Russian and with a preface by E. R. Dawson,
  Wiley-Interscience Series in Discrete Mathematics.

\bibitem{PTZ12}
Richard Peng, Kanat Tangwongsan, and Peng Zhang.
\newblock Faster and simpler width-independent parallel algorithms for positive
  semidefinite programming, 2012.

\bibitem{small1990survey}
Christopher~G Small.
\newblock A survey of multidimensional medians.
\newblock {\em International Statistical Review/Revue Internationale de
  Statistique}, pages 263--277, 1990.

\bibitem{MR0120720}
John~W. Tukey.
\newblock A survey of sampling from contaminated distributions.
\newblock In {\em Contributions to probability and statistics}, pages 448--485.
  Stanford Univ. Press, Stanford, Calif., 1960.

\bibitem{MR0133937}
John~W. Tukey.
\newblock The future of data analysis.
\newblock {\em Ann. Math. Statist.}, 33:1--67, 1962.

\end{thebibliography}
\end{footnotesize}

\end{document}